\documentclass[final,onefignum,onetabnum]{siamart171218}

\usepackage{amsfonts,graphicx}
\usepackage[caption=false]{subfig}
\usepackage{epstopdf,algorithmic}
\usepackage{amsmath}
\DeclareMathOperator*{\argmax}{argmax}
\usepackage{tikz,pgfplots}
\usepackage{booktabs}
\usetikzlibrary{patterns}
\usepgfplotslibrary{groupplots}
\pgfplotsset{compat=1.17}
\pgfplotsset{
    every non boxed x axis/.style={} 
}

\definecolor{hcorange}{RGB}{245, 130, 48}
\definecolor{hcnavy}{RGB}{0, 0, 128}
\definecolor{hcblue}{RGB}{0, 130, 200}
\definecolor{hcpink}{RGB}{250, 190, 190}
\definecolor{hcbrown}{RGB}{128, 0, 0}
\definecolor{hclavender}{RGB}{230, 190, 255}
\definecolor{hcgrey}{RGB}{128, 128, 128}
\definecolor{hcgreen}{RGB}{60, 180, 75}
\definecolor{hcred}{RGB}{230, 25, 75}

\definecolor{red1_1}{RGB}{162, 13, 30}
\definecolor{red2_1}{RGB}{188, 13, 53}
\definecolor{red3_1}{RGB}{222, 38, 76}
\definecolor{red4_1}{RGB}{240, 120, 140}
\definecolor{red5_1}{RGB}{246, 177, 197}
\definecolor{green1_2}{RGB}{0, 57, 66}
\definecolor{green2_2}{RGB}{0, 90, 92}
\definecolor{green3_2}{RGB}{0, 117, 108}
\definecolor{green4_2}{RGB}{0, 138, 112}
\definecolor{green5_2}{RGB}{2, 167, 117}
\definecolor{blue1_2}{RGB}{0, 48, 90}
\definecolor{blue2_2}{RGB}{0, 75, 141}
\definecolor{blue3_2}{RGB}{0, 117, 219}
\definecolor{blue4_2}{RGB}{65, 146, 217}
\definecolor{blue5_2}{RGB}{122, 186, 242}

\usepackage{longtable,tabularx}
\newcolumntype{M}[1]{>{\centering\arraybackslash}m{#1}}
\usepackage{threeparttable}
\ifpdf
  \DeclareGraphicsExtensions{.eps,.pdf,.png,.jpg}
\else
  \DeclareGraphicsExtensions{.eps}
\fi

\newsiamremark{remark}{Remark}
\newsiamremark{hypothesis}{Hypothesis}
\crefname{hypothesis}{Hypothesis}{Hypotheses}
\newsiamthm{claim}{Claim}

\headers{Numerically Stable s-step GMRES}{Z. Xu, J. J. Alonso, and E. Darve}

\title{A numerically stable communication avoiding s-step GMRES algorithm\thanks{Submitted to the editors \today.}
\funding{This work was funded by NASA Grant 80NSSC18M0152 from the NASA Transformational Tools and Technologies program.}}

\author{Zan Xu\thanks{Department of Aeronautics and Astronautics, Stanford University, Stanford, CA 94305
  (\email{zanxu@stanford.edu}, \email{jjalonso@stanford.edu}).}
\and Juan J.~Alonso\footnotemark[2]
\and Eric Darve\thanks{Department of Mechanical Engineering, Stanford University, Stanford, CA 94305
  (\email{darve@stanford.edu}).}}

\usepackage{amsopn}

\makeatletter
\newcommand*{\addFileDependency}[1]{% argument=file name and extension
  \typeout{(#1)}% latexmk will find this if $recorder=0 (however, in that case, it will ignore #1 if it is a .aux or .pdf file etc and it exists! if it doesn't exist, it will appear in the list of dependents regardless)
  \@addtofilelist{#1}% if you want it to appear in \listfiles, not really necessary and latexmk doesn't use this
  \IfFileExists{#1}{}{\typeout{No file #1.}}% latexmk will find this message if #1 doesn't exist (yet)
}
\makeatother

\usepackage{cleveref}
\crefname{equation}{}{}
\Crefname{equation}{}{}
\crefname{figure}{fig.}{figs.}
\Crefname{figure}{Fig.}{Figs.}
\crefname{algorithm}{algorithm}{algorithms}
\Crefname{algorithm}{Algorithm}{Algorithms}
\crefname{subsubsection}{test}{tests}
\Crefname{subsubsection}{Test}{Tests}

\usepackage{chngcntr}

\counterwithout{subsubsection}{subsection} % remove section from the subsidiary counter

\usepackage{titlesec}
\titleformat{\subsubsection}[runin]{\normalfont\normalsize\bfseries}{Test~\thesubsubsection: }{0em}{}
\PassOptionsToPackage{normalem}{ulem}
\usepackage{ulem}

\newcommand{\rpy}[1]{{\color{red}{}#1}} % reply

\usetikzlibrary{external}

\ifpdf
\hypersetup{
  pdftitle={A numerically stable communication-avoiding \texorpdfstring{$s$}{s}-step GMRES algorithm},
  pdfauthor={Zan Xu, Juan J.~Alonso, and Eric Darve}
}
\fi

\begin{document}

% first round of review
% \input{Review/cover_letter.tex} % this contains the cover letter to editor
% \input{Review/review.tex} % this contains response to all reviewers / editor
% \input{Review/reviewer1.tex} % this contains response specifically to reviewer 1
% \input{Review/reviewer2.tex} % this contains response specifically to reviewer 2

% second round of review
%\input{Review/reviewer1_round2.tex} % this contains response specifically to reviewer 2
%\input{Review/reviewer2_round2.tex} % this contains response specifically to reviewer 2
%
%% % reset a few things after the review response.
%\newpage
%\setcounter{page}{1}
%\resetlinenumber

\maketitle

% REQUIRED
\begin{abstract}
  Krylov subspace methods are extensively used in scientific computing to solve large-scale linear systems. However, the performance of these iterative Krylov solvers on modern supercomputers is limited by expensive communication costs. The $s$-step strategy generates a series of $s$ Krylov vectors at a time to avoid communication. Asymptotically, the $s$-step approach can reduce communication latency by a factor of $s$. Unfortunately, due to finite-precision implementation, the step size has to be kept small for stability. In this work, we tackle the numerical instabilities encountered in the $s$-step GMRES algorithm. By choosing an appropriate polynomial basis and block orthogonalization schemes, we construct a communication avoiding $s$-step GMRES algorithm that automatically selects the optimal step size to ensure numerical stability. 
  To further maximize communication savings, we introduce scaled Newton polynomials that can increase the step size $s$ to a few hundreds for many problems. An initial step size estimator is also developed to efficiently choose the optimal step size for stability. The guaranteed stability of the proposed algorithm is demonstrated using numerical experiments. In the process, we also evaluate how the choice of polynomial and preconditioning affects the stability limit of the algorithm. Finally, we show parallel scalability on more than 114,000 cores in a distributed-memory setting. Perfectly linear scaling has been observed in both strong and weak scaling studies with negligible communication costs.
\end{abstract}

% REQUIRED

% REQUIRED

\section{Introduction}
Krylov subspace methods are widely used to find iterative solutions to linear systems of equations. They are particularly attractive for sparse and large problems, where direct methods are often prohibitively expensive in terms of both computation and storage. However, on modern large-scale computing architectures, the performance of iterative Krylov solvers is limited by expensive communication costs~\cite{demmel2008avoiding, mohiyuddin2009minimizing}. Such costs are incurred when moving data either between levels of the local memory hierarchy in a single processor or between different processors connected over a network. The performance of many operations in Krylov subspace methods, such as dot products that require global reductions, is bounded by communication on distributed systems~\cite{hoemmen2010communication}. As technology trends project floating-point performance to continue to scale at a faster rate than latency reduction~\cite{national2005getting}, communication costs in numerical algorithms are becoming a serious bottleneck. Hence, to effectively leverage new high-performance computing hardware, one needs to redesign Krylov subspace methods to communicate less.

Communication-avoiding $s$-step Krylov subspace methods~\cite{bai1994newton, chronopoulos1989s,chronopoulos1991s,chronopoulos1996parallel,de1995reducing,hoemmen2010communication, joubert1992parallelizable, carson2015communication, carson2018adaptive, mohiyuddin2009minimizing} attempt to \textit{avoid} communication by generating a series of $s$ Krylov basis vectors at once and orthogonalize them in a block fashion. This strategy replaces most BLAS-1 and BLAS-2 operations with BLAS-3 operations that are better optimized for cache memory~\cite{dongarra1990set}. Furthermore, the appropriate choice of block orthogonalization schemes can reduce the number of global communication constructs by a factor of $s$~\cite{hoemmen2010communication}. On parallel computers, this not only decreases global communication latency by $s$ asymptotically but also reduces the number of synchronization points in the algorithm, thus mitigating system noise, load imbalances and/or hardware variations among different processors~\cite{swirydowicz2020low,yamazaki2020low}.

$S$-step Krylov subspace methods are mathematically equivalent to their classical counterparts. However, finite-precision implementations of these variants introduce several sources of numerical instability, such as the poor conditioning of the Krylov basis matrix~\cite{beckermann2000condition,philippe2012generation} and the loss of orthogonality of basis vectors due to block QR schemes~\cite{carson2015communication,carson2022block}. As a result, the maximal step size $s$ in $s$-step methods has to be limited to small values for convergence. Since communication savings rely directly on the step size, it is necessary to address stability issues to improve the scalability of these methods on high-performance computing resources.

This paper focuses on $s$-step variants of the Generalized Minimal Residual (GMRES) algorithm~\cite{saad1986gmres}. We first review the stability and communication pattern of $s$-step GMRES algorithms in \Cref{sec:background}, and highlight the limitations in the current state of the art. We then derive an adaptive $s$-step GMRES algorithm in \Cref{sec:algorithm}, which dynamically determines the maximum allowed step size $s$ in a communication-avoiding manner. To achieve better communication savings, we also introduce 
scaled Newton polynomials in \Cref{subsec:scalednewton} that can stably increase the step size $s$ to a few hundred for many problems, and an initial step size estimator in \Cref{subsec:errorestimator} to efficiently choose initial step sizes. \Cref{sec:tests_stability} presents a series of numerical experiments that assess the numerical stability of the proposed algorithm under different measures. Finally, we demonstrate the parallel scalability of the overall algorithm on more than 114,000 cores in \Cref{sec:scalability} before we conclude with some remarks in \Cref{sec:conclusions}.

\section{Background on \texorpdfstring{$s$}{s}-step GMRES}
\label{sec:background}

Same as the conventional GMRES algorithm, the $s$-step GMRES method solves the linear system
\begin{equation*}
    Ax = b, 
\end{equation*}
with $A \in \mathbb{R}^{N\times N}$. At each iteration, the Matrix Powers Kernel (MPK) performs $s$ SpMVs on a starting normalized column vector $Q_i$ at block iteration $i$ to generate a Krylov basis matrix $V(Q_i,s)$ containing $s+1$ Krylov basis vectors
\begin{equation*}
    V(Q_i, s) = \big[p_0(A)Q_i, \, p_1(A)Q_i, \, p_2(A)Q_i, \, p_3(A)Q_i,\, \dots \,,\, p_s(A)Q_i\big] \in \mathbb{R}^{N\times (s+1)},
\end{equation*}
where $s$ is the step size (or block size\footnote{Throughout this work, we use \textit{block size} and \textit{step size} interchangeably to refer to multiple Krylov basis vectors that are generated at the same time. This terminology does not have the same connotation as \textit{block size} in block-GMRES literature, where each block refers to multiple right-hand sides.}) at each block iteration, and $p_j(z)$ is a polynomial of degree $j$ defined by a three-term recurrence
\begin{align*}
\begin{split}
    p_0(z) &= 1, \quad p_1(z) = (z-\theta_0)p_0(z)/\gamma_0, \\
    p_j(z) &= ((z-\theta_{j-1})p_{j-1}(z) - \beta_{j-2} p_{j-2}(z))/\gamma_{j-1}.
\end{split}
\end{align*}
Parameters of the polynomial function $(\theta_i, \beta_i, \gamma_i)$ can be condensed in a \textit{change of basis} matrix $B$, which is used as input to the MPK~\cite{hoemmen2010communication,ballard2014communication,carson2015communication}. Communication-avoiding MPK can be achieved~\cite{demmel2008avoiding,hoemmen2010communication} though standard SpMVs are often used to accommodate different preconditioners. The Krylov basis matrix $V(Q_i,s)$ is then orthogonalized using an appropriate QR factorization that ought to require as little communication as possible. The remaining downstream operations related to the upper Hessenberg matrix $H$ and the least-square problem are local to each processor, thus not affecting the communication pattern of the algorithm. The overall framework of the $s$-step GMRES is presented in \Cref{alg:sgmres}, and detailed algorithm implementations can be found in~\cite{hoemmen2010communication,ballard2014communication,yamazaki2017improving}.

\begin{algorithm}[ht]
\caption{$s$-step GMRES}\label{alg:sgmres}
\begin{algorithmic}[1]
\REQUIRE $N\times N$ matrix $A$, right-hand side vector $b$, initial guess vector $x_0$, maximum iteration count $m$, step size \textit{s}, change of basis matrix $B$.
\ENSURE $x$, solution to the linear system $Ax=b$.
\STATE $r := b - Ax_0, \quad \beta := \|r\|_2, \quad Q_1 := r/\beta, \quad R_{1,1}=1, \quad i = 1$
\WHILE{$i \leq m$}
\STATE $V$ = MPK$(A, Q_i, s, B)$
\STATE $[Q_{i+1:i+s}, R_{1:i+s,i+1:i+s}]$ = Block QR($V, Q_{1:i}, s$)
\STATE Assemble upper Hessenberg matrix $H_{1:i+s,1:i+s-1}$
\STATE Apply Givens rotation to update matrix $H_{1:i+s,1:i+s-1}$
\STATE Check for convergence
\STATE $i = i + s$
\ENDWHILE
\STATE$y = $ argmin$\|\beta e_1 - H_{1:m+1,1:m}y\|_2$ \COMMENT{Least-squares problem}
\STATE$x = x_0 + Q_{1:m} y$
\end{algorithmic}
\end{algorithm}

\subsection{Numerical instability in the MPK}
\label{subsec:stability_MPK}
The first source of numerical instability of the $s$-step GMRES algorithm lies in the poor conditioning of the Krylov basis matrix $V(Q_i,s)$ generated from the MPK. The early work on $s$-step methods~\cite{walker1988implementation, chronopoulos1991s} started with the monomial basis,
\begin{equation}\label{eq:monomial}
    p_j(A) = A^j.
\end{equation}
But in finite precision, the Krylov basis vectors converge to the eigenvector associated with the largest eigenvalue of matrix $A$, making $V(Q_i,s)$ ill-conditioned or even numerically rank-deficient. Moreover, to avoid overflow or underflow, each vector obtained from SpMV needs to be normalized, leading to $s$ additional global communications that are undesirable. 

To slow the growth of the incremental condition number of $V(Q_i,s)$ for large values of $s$, many approaches have been tried. Bai et al.\ proposed a Newton basis GMRES implementation using approximate eigenvalues of $A$ as $\theta_i$'s~\cite{bai1994newton,hoemmen2010communication}. The Chebyshev basis has also been studied and has been shown to slow the growth of condition numbers with varying degrees of success~\cite{joubert1992parallelizable,joubert1992parallelizable2,philippe2012generation,ballard2014communication,carson2015communication}. The normalization step of the basis vectors can be eliminated in order to reduce global communications if $s$ is kept relatively small. Without normalization, methods such as matrix equilibration have been adopted to slow the growth/decay of Krylov basis vector norms and increase $s$~\cite{hoemmen2010communication,ballard2014communication}. 

However, since most approaches are heuristic in nature, the maximum allowed value of $s$ for a given input matrix $A$ can only be determined by trial and error or operations involving expensive communication constructs. Because of the lack of an appropriate mechanism to monitor the conditioning of $V(Q_i,s)$, the ill-conditioning of the Krylov basis matrix $V(Q_i,s)$ poses challenges to the downstream QR factorization and undermines the stability of the overall algorithm.

\subsection{Numerical instability in block QR orthogonalization}
\label{subsec:stability_QR}
Block QR me\-th\-ods are commonly used to orthonormalize Krylov basis vectors in $V(Q_i,s)$ at the same time, giving the potential to reduce communication latency by a factor of $s$. To set the stage for a detailed discussion, we first introduce the formal notation for block QR orthogonalization.

For a matrix $\mathcal{V} \in \mathbb{R}^{N \times M}$, partitioned into $p$ blocks,
\begin{equation}\label{eq:matrixV}
    \mathcal{V} = [\,V^{(1)}\,|\,V^{(2)}\,|\, \dots \,|\,V^{(p)}\,],
\end{equation}
we seek an ``economic" block QR factorization in the form of 
\begin{equation}
    \mathcal{V} = \mathcal{Q} R,
\end{equation}
with $\mathcal{Q} = [\,Q^{(1)}\,|\,Q^{(2)}\,|\, \dots \,|\,Q^{(p)}\,] \in \mathbb{R}^{N \times M}$ having the same block structure as $\mathcal{V}$ and $R\in\mathbb{R}^{M \times M}$. 
In this work, we focus on the \textit{loss of orthogonality} (LOO) error defined by  
\begin{equation}
    \|I - \mathcal{Q}^T\mathcal{Q}\|_F,
\end{equation}
as a measure of stability for QR factorization. Typically, block QR consists of two stages in each block iteration. Each block $\mathcal{V}^{(i)}$ first needs to be \textit{inter-orthogonalized} with respect to orthogonal bases generated from previous blocks ($Q^{(1)}$ to $Q^{(i-1)}$), then \textit{intra-orthogonalized} to ensure orthogonality among vectors within the same block.

The current state-of-the-art communication-avoiding $s$-step GMRES focuses on block classical Gram-Schmidt with reorthogonalization (BCGS2) for inter-ortho\-go\-na\-li\-zation and Cholesky QR (CholQR\footnote{Due to re-orthogonalization in inter-orthogonalization, CholQR is applied twice in total. However, we use the abbreviation CholQR instead of CholQR2 as each instance of intra-orthogonalization involves applying CholQR only once. We reserve the term CholQR2 for referring to applying CholQR twice \textit{within} one instance of intra-ortho\-go\-na\-li\-zation.}) for intra-orthogonalization (\Cref{alg:BlockQR}). Block orthogonalization routines often utilize cache-friendly BLAS-3 kernels to minimize data movement among memory hierarchies~\cite{golub2013matrix}. In a distributed-memory setting, BCGS and CholQR only cost one global synchronization each in one block iteration. However, a re-orthogonalization step is often necessary to reduce the LOO error~\cite{barlow2013reorthogonalized,yamamoto2015roundoff}, doubling the cost of internode communication to four global synchronizations per block iteration. In~\cite{yamazaki2020low} a low-synchronization variant of BCGS2 with CholQR has been proposed that uses one global reduction per $s$ vectors. 

\begin{algorithm}
\caption{Block QR in \Cref{alg:sgmres} using BCGS2 and CholQR}  \label{alg:BlockQR}
\begin{algorithmic}[1]
\REQUIRE Step size \textit{s}, index $i$, Krylov basis matrix $V$ of size $N\times s$, orthogonal basis matrix $Q_{1:i}$ of size $N\times i$.
\ENSURE Orthogonal basis matrix $Q_{i+1:i+s}$ of size $N\times s$, upper triangular matrix $R_{1:i+s,i+1:i+s}$ of size $(i+1) \times s$.
\STATE $k = i+1:i+s$
\STATE
\STATE $W = Q_{1:i}^TV$ \COMMENT{First BCGS}
\STATE $V = V- Q_{1:i}W$
\STATE $Z = \mbox{Cholesky}(V^TV)$ \COMMENT{First CholQR}
\STATE $Q_{k} = V/Z$
\STATE
\STATE $R_{1:i, k} = Q_{1:i}^TQ_{k}$ \COMMENT{Second BCGS}
\STATE $Q_{k} = Q_{k} - Q_{1:i}R_{1:i, k}$
\STATE $\tilde{Z} = \mbox{Cholesky}(Q_{k}^TQ_{k})$ \COMMENT{Second CholQR}
\STATE $Q_{k} = Q_{k}/\tilde{Z}$
\STATE
\STATE $R_{1:i, k} = W + R_{1:i, k}Z$ \COMMENT{Combine both steps}
\STATE $R_{k, k} = \tilde{Z}Z$
\end{algorithmic}
\end{algorithm}

The reduction in communication cost of block QR methods has made them attractive for large-scale problems. However, for most block QR methods, formal proofs of their stability properties remain elusive. We refer the reader to~\cite{carson2022block} for a recent survey on the stability properties of common block QR methods based on Gram-Schmidt algorithms. The numerical experiments therein (e.g., Figs.\ 5 and 6) demonstrated that even the current state of the art is not unconditionally stable for handling arbitrary step sizes with finite precision. In this work, we focus primarily on the stability analysis of the original BCGS2 with the CholQR algorithm. The stability analysis of low-synchronization variants is left for future work.

For a numerically non-singular matrix $\mathcal{V}$ with condition number bound $\kappa(\mathcal{V}) < O(\epsilon^{-1})$, where $\epsilon$ is the unit round-off, Barlow showed that BCGS2 produces $O(\epsilon)$ orthogonality error if the intra-orthogonalization scheme can achieve the same level of orthogonality~\cite{barlow2013reorthogonalized}. The round-off error analysis performed in~\cite{yamamoto2015roundoff} indicates that CholQR applied twice does give the desired $O(\epsilon)$ level of LOO error but requires the input block to have a condition number below $O(\epsilon^{-1/2})$. Here, we note that the input block does not refer to the block $V^{(i)}$ in the original matrix $\mathcal{V}$ in~\Cref{eq:matrixV}, but to the block matrix used to generate the Gram matrix. (For example, notice how line 4 updates the block matrix $V$ in \Cref{alg:BlockQR} before the CholQR step.) Without any safeguard on the conditioning of the input matrix, the output of CholQR can potentially degrade the numerical stability and lead to a significant increase in the LOO error of the overall block QR algorithm. Furthermore, we emphasize that even if efforts were made to improve and constrain the condition number of the Krylov basis matrix $V(Q_i,s)$ in the upstream MPK, the inter-orthogonalization process could nonetheless lead to an ill-conditioned input for the intra-orthogonalization. 

Sometimes, in a more catastrophic setting, an ill-conditioned input block matrix in CholQR (e.g., $V$ in line 5 of~\Cref{alg:BlockQR}) can give rise to a Gram matrix that is not symmetric positive-definite (SPD) in finite precision, which will cause the Cholesky decomposition to break down. Such runtime errors often force the user to switch to other QR methods that are not communication-avoiding or restart the algorithm with a smaller block size $s$. 

\subsection{Limitations of current state-of-the-art methods} \label{subsec:limitations} As the communication savings of $s$-step GMRES are directly related to the step size $s$, one would like to use a step size as large as possible. However, due to the numerical instabilities mentioned above, conservative values of $s$ (e.g., $5$, $10$) are often reported in the literature~\cite{hoemmen2010communication,ballard2014communication,yamazaki2017improving,yamazaki2020low} to mitigate numerical problems. On large-scale distributed computing architectures, finding the maximum allowed choice of $s$ of an arbitrary matrix by trial and error is highly inefficient and costly. To the best of our knowledge, no variant of $s$-step GMRES in the existing literature can determine a stable value of $s$ \textit{a priori} in a communication-avoiding fashion, nor can $s$ be dynamically adapted to ensure stability. Attempts were made to improve the stability of $s$-step GMRES using a predetermined sequence of block sizes, which is difficult to generalize~\cite{imberti2017varying}. An adaptive mechanism has been derived for the $s$-step conjugate gradient (CG) algorithm~\cite{carson2018adaptive} but does not translate to the $s$-step GMRES easily. One key difference is that the $s$-step CG method does not require an explicit block QR algorithm. Furthermore, the adaptive mechanism in $s$-step CG relies on the norm of the residual vectors, which is only available in $s$-step GMRES \textit{after} the block QR is completed. Thus, the numerical instabilities described in \Cref{subsec:stability_QR} cannot be resolved using the same approach. By and large, the backward stability of $s$-step GMRES still remains an open problem, and the numerical instabilities greatly limit the potential communication savings in GMRES algorithms.

\section{Proposed solutions to resolve numerical instabilities in \texorpdfstring{$s$}{s}-step GMRES}
\label{sec:algorithm}
In this section, we aim to address the numerical instability in the current state-of-the-art $s$-step GMRES and derive an adaptive $s$-step GMRES algorithm that is numerically stable by construction. In particular, we target the $s$-step GMRES algorithm that uses BCGS2 for inter-orthogonalization and CholQR for intra-orthogonalization. We start by resolving the breakdown problem in CholQR in \Cref{subsec:pchol}, which provides a path for us to design a dynamical mechanism in \Cref{subsec:stopping} that monitors the step size $s$ for stability. We then provide the complete algorithm in \Cref{subsec:adaptsgmres} and discuss its communication patterns and computational trade-offs. To further improve communication savings, we introduce \textit{scaled} Newton polynomials in~\Cref{subsec:scalednewton} as a basis to increase the step size $s$ and an initial step size estimator in~\Cref{subsec:errorestimator} to pick the optimal step size at the first iteration.

\subsection{Avoiding Cholesky breakdown in CholQR}\label{subsec:pchol}

As mentioned previously, an ill-conditioned input matrix to CholQR (e.g., $V$ in line 5 of~\Cref{alg:BlockQR}) may produce a numerically non-SPD Gram matrix that causes the Cholesky decomposition to break down. To avoid this, we propose the use of a \textit{partial} Cholesky decomposition, which leads to a \textit{partial} CholQR factorization. 

Let $X \in \mathbb{R}^{N\times s}$ be the input to the CholQR factorization, and let $p$ be an integer such that
\begin{equation}
    X = [X_p\, |\, X_{s-p}],\qquad 1\leq p \leq s,
\end{equation}
where $X_p \in \mathbb{R}^{N\times p}$ represents the first $p$ columns of $X$ and $X_{s-p} \in \mathbb{R}^{N\times (s-p)}$ the remaining columns.\footnote{For brevity, the subscript in this section represents the number of rows/columns. This is to be differentiated from the previous notation, where a single subscript denotes a column vector (e.g., $Q_i$).} The first step of CholQR computes a Gram matrix $X^TX$,
\begin{equation}
    X^TX = \begin{bmatrix}
    X_p^TX_p & X_p^TX_{s-p} \\
    X_{s-p}^TX_p & X_{s-p}^TX_{s-p} \\
    \end{bmatrix}.
\end{equation}
If we represent $R$ as the upper triangular matrix obtained from Cholesky decomposition of $X^TX$,
\begin{equation}
    R = \text{Cholesky}(X^TX) = \begin{bmatrix}
    R_p & R_{p,s-p} \\
    0 & R_{s-p} \\
    \end{bmatrix},
\end{equation}
it follows that $R_p = \text{Cholesky}(X_p^TX_p)$. Hence, by using the triangular matrix $R_p$ obtained from Cholesky decomposition of $X_p^TX_p$ ($p^{th}$ order leading principal submatrix of $X^TX$), one can compute a \textit{partial} CholQR factorization of the matrix $X$ to obtain $p$ columns of the orthogonal basis $Q_p$:
\begin{equation}\label{eq:pchol}
\begin{split}
    Q_p & = X_p R_p^{-1}, \\
    X \, &= [Q_p R_p, X_{s-p}]. \\
\end{split}
\end{equation}
This formulation gives flexibility to stop the QR factorization at an arbitrary column of the input matrix $X$. In finite precision, this strategy can be used to avoid a Cholesky breakdown, as one can always salvage the factorization computed right before the breakdown. In the worst-case scenario of $p=1$ (i.e., Cholesky breaks down at the second row/column), partial CholQR is reduced to a norm computation, which is a backward stable floating point operation~\cite{higham2002accuracy}. 

The remaining columns, $X_{s-p}$, can be recursively orthogonalized. However, in the context of $s$-step GMRES, the remaining columns correspond to higher powers of a polynomial basis and are often harder to orthogonalize, leading to small values of block size $p$ in subsequent iterations. Sometimes, due to the ill-conditioned Krylov matrix from upstream MPK, the remaining columns can be linearly dependent on the previous block $X_p$ and are thus not worth orthogonalizing as they do not contribute to the construction of the Krylov subspace. Therefore, to maximize communication savings, the remaining columns are discarded, and the last column of $Q_p$ can be used to generate the next Krylov basis matrix in the next block iteration. This is justified on large-scale distributed machines as we avoid inefficient Cholesky breakdowns at the expense of some SpMVs that are cheaper in terms of communication costs.

\subsection{Stopping criteria for partial CholQR} \label{subsec:stopping}
A standard implementation of Cholesky decomposition breaks down when the diagonal entry to be factorized is equal to or less than zero. The partial CholQR introduced in \Cref{subsec:pchol} is sufficient to avoid Cholesky breakdowns by salvaging the calculated factorization. However, this does not guarantee minimal LOO. In other words, stopping at non-positive pivots in Cholesky is insufficient to fulfill the condition number bound for CholQR. A more stringent choice of $p$ in \Cref{eq:pchol} is necessary to ensure that the columns $X_p$ have a condition number at most $O(\epsilon^{-1/2})$,
\begin{equation}\label{eq:pchol2}
    \kappa(X_p) = \kappa(R_p) \leq O(\epsilon^{-1/2}).
\end{equation}
As a result, an appropriate stopping criterion must be imposed after each rank-1 Cholesky downdate to bound the condition number growth while generating $R_p$. 

Expensive singular value decomposition (SVD) can be used to compute condition numbers for small dense matrices and was suggested in the rank-revealing TSQR (Tall-skinny QR) algorithm~\cite{hoemmen2010communication} for condition number monitoring. However, since the algorithmic complexity to compute the SVD of $R_p$ is $O(p^3)$ and the stopping criterion must be invoked after each rank-1 update, the total cost to check against the stopping criterion is on the order of $O(p^4)$, exceeding the $O(p^3)$ complexity of Cholesky decomposition. This could be justified for small block sizes. 

If one wants to use large block sizes, the incremental condition estimator (ICE)~\cite{bischof1990incremental,bischof1991robust} can be used as a cheaper alternative. After each rank-1 update, ICE updates its condition number estimation using the newly computed column of $R_p$ without accessing the entire leading principal submatrix. Consequently, ICE costs only $O(p)$ flops at each iteration. An example of the resulting partial Cholesky decomposition is illustrated in \Cref{alg:pcholICE}. 

We recognize that the condition number estimate provided by ICE is often an underestimation, but the discrepancy is generally within a factor of $10$~\cite{bischof1991robust}. In practice, we found that a condition number upper bound of $\Omega \sim O(10^{-1}\epsilon^{-1/2})$ is sufficient for \Cref{eq:pchol2} to hold numerically.

\begin{algorithm}
\caption{Partial Cholesky decomposition}  \label{alg:pcholICE}
\begin{algorithmic}[1]
\REQUIRE $m\times m$ matrix $A$, and condition number upper bound $\Omega$
\ENSURE $j$, and $j \times j$ triangular matrix $R$ such that $A_{1:j,1:j} = R^TR$
\STATE $\kappa_0 = 1$
\FOR[Standard Cholesky]{$j = 1, ..., m$}
\FOR{$i = 1,...,j-1$}
\STATE$R_{i,j} = (A_{i,j} - \sum_{k=1}^{j-1}R_{k,i}R_{k,j})/R_{i,i}$
\ENDFOR
\STATE$R_{j,j} = \sqrt{A_{j,j} - \sum_{k=1}^{j-1}R_{k,j}^2}$
\IF{SVD}
\STATE $\sigma_{\max}, \sigma_{\min} = {\rm SVD}(R_{1:j,1:j})$ \COMMENT{Singular Value Decomposition}
\STATE $\kappa_j = \sigma_{\max}/\sigma_{\min}$
\ELSIF{ICE} 
\STATE $\kappa_j = {\rm ICE}(R_{1:j,j}, \kappa_{j-1})$  \COMMENT{Incremental Condition Estimator}
\ENDIF
\IF{$\kappa_j > \Omega$}
\STATE $j = j - 1$ 
\STATE break
\ENDIF
\ENDFOR
\end{algorithmic}
\end{algorithm}

\subsection{Adaptive \texorpdfstring{$s$}{s}-step GMRES algorithm}\label{subsec:adaptsgmres}

Using BCGS2 with a partial Chol\-QR factorization in $s$-step GMRES, we can dynamically monitor the growth of the condition number and adjust the block size $s$ to ensure that the orthogonalization process stays within the bounds of the condition number. An outline of the resulting adaptive $s$-step GMRES is presented in \Cref{alg:adaptsgmres}. We make a few remarks on the overall algorithm's communication pattern, stability, and other practical aspects.

\begin{algorithm}
\caption{Adaptive $s$-step GMRES}  \label{alg:adaptsgmres}
\begin{algorithmic}[1]
\REQUIRE $N\times N$ matrix $A$, right hand side vector $b$, initial guess vector $x_0$, maximum iteration count $m$, \textit{initial} step size \textit{$s_0$}, change of basis matrix $B$, condition number upper bound $\Omega$.
\ENSURE $x$, solution to the linear system $Ax=b$.
\STATE $r := b - Ax_0, \quad Q_{1} := r/\|r\|_2, \quad R_{1,1} = 1,  \quad i = 1, \quad s = s_0$
\WHILE{$i \leq m$}
\STATE $V$ = MPK($A, Q_{i}, B, s$) 
\STATE
\STATE $W = Q_{1:i}^TV$ \COMMENT{First BCGS}
\STATE $V = V- Q_{1:i}W$
\STATE $[\,p,\, Z\,] = \mbox{Partial Cholesky}(V^TV, \Omega)$ \COMMENT{First partial CholQR}
\STATE $k = i+1:i+p$
\STATE $Q_{k} = V_{:,1:p}/Z_{1:p,1:p}$
\STATE
\STATE $R_{1:i, k} = Q_{1:i}^TQ_{k}$ \COMMENT{Second BCGS}
\STATE $Q_{k} = Q_{k} - Q_{1:i}R_{1:i, k}$
\STATE $[\, p\,, \tilde{Z}\,] = \mbox{Partial Cholesky}(Q_{k}^TQ_{k}, \Omega)$ \COMMENT{Second partial CholQR}
\STATE $k = i+1:i+p$
\STATE $Q_{k} = Q_{k}/\tilde{Z}$
\STATE
\STATE $R_{1:i, k} = W_{:,1:p} + R_{1:i, k}Z_{1:p,1:p}$ \COMMENT{Combine both steps}
\STATE $R_{k, k} = \tilde{Z}Z_{1:p,1:p}$
\STATE
\STATE $s = p$ \COMMENT{Update step size}
\STATE
\STATE Assemble upper Hessenberg matrix $H_{1:i+s,i:i+s-1}$
\STATE Apply Givens rotation to update matrix $H_{1:i+s,i:i+s-1}$
\STATE Check for convergence
\STATE $i = i + s$
\ENDWHILE
\STATE$y = $ argmin$\|\beta e_1 - H_{1:m+1,1:m}y\|_2$ \COMMENT{Least-squares problem}
\STATE$x = x_0 + Q_{1:m} y$
\end{algorithmic}
\end{algorithm}

The formal proof of the backward stability of the overall algorithm is reserved for future work. Here, we provide some conjectures on stability based on established results and verify the stability using numerical experiments in \Cref{sec:tests_stability}.

The discussion in Section 4.2 of~\cite{carson2022block}, which builds on the analysis in~\cite{barlow2013reorthogonalized}, has already indicated that BCGS2 with a merely conditionally stable intra-orthogonalization scheme such as Classical Gram-Schmidt (CGS) can still give $O(\epsilon)$ LOO error for well-conditioned matrices. Since CholQR has the same stability bound and orthogonality error behavior as CGS, we conjecture that BCGS2 with CholQR is stable for well-conditioned matrices. For challenging matrices, authors in~\cite{carson2022block} concluded that unconditionally stable intra-orthogonalization schemes, such as Householder QR, are necessary for the overall stability, under the key assumption that each block to be factorized be of a fixed size $s$. In our work, using partial CholQR, we relax such an assumption and allow the algorithm to adjust the step size on the fly. This has several advantages from a numerical viewpoint.
\begin{itemize}
    \item As the block iteration continues, the subsequent blocks generally tend to be harder to orthogonalize than the first few blocks, as the inter-orthogonali\-za\-tion step is looking for new sets of orthogonal bases within the remaining subspace, which generally gets smaller with more block iterations. In such cases, the adaptive mechanism allows smaller step sizes for later blocks.
    \item For very ill-conditioned matrices where the adapted step size turns out to be $s = 1$, the partial CholQR is reduced to a norm computation. The overall block QR factorization is reduced to a column-wise CGS2 (Classical Gram-Schmidt with re-orthogonalization) algorithm, which is known to be stable either as a standalone QR factorization~\cite{daniel1976reorthogonalization,drkovsova1995numerical} or as an orthogonalization scheme used in GMRES~\cite{swirydowicz2020low}.
    \item With a large block size $s$, the MPK can generate Krylov matrices $V(Q_i, s)$ that could be numerically rank-deficient. Under such circumstances, there is no viable QR factorization, including unconditionally stable QR factorization schemes such as Householder QR that can expand the dimension of the Krylov subspace by $s$. The ability to reduce the block size $s$ and discard the remaining columns corresponding to higher powers of polynomial basis in partial CholQR avoids such problems. The algorithm always works on the maximum allowed set of Krylov vectors for orthogonalization.
\end{itemize}

By discarding some basis vectors in the Krylov basis matrix $V(Q_i, s)$, we acknowledge that some SpMV computations done earlier in the MPK are potentially wasted. The same issue was also reported in the adaptive $s$-step conjugate gradient method~\cite{carson2018adaptive}. Here, we argue that this is still desirable in a large-scale distributed-memory system for many applications (e.g., solutions to PDEs) where SpMVs require mostly Point-to-Point communication with neighboring nodes, but global synchronizations involve much more expensive collective communication costs. To minimize wasted effort, we update the step size at each block iteration (line 20 in \Cref{alg:adaptsgmres}) since the maximum allowable step size is generally non-increasing. The initial step size, however, at the very first iteration is a user-defined value. Poor choice of the initial step size can still result in inefficient use of the algorithm. A step size that is too small may not maximize communication savings, while a step size that is too large may result in truncation of many Krylov basis vectors, leading to a significant amount of wasted computations. 

To this end, the following subsections look into ways to maximize communication savings and minimizing wasted computations. We first introduce a polynomial basis in \Cref{subsec:scalednewton} that allows for very large step sizes in many problems and then derive an estimator in \Cref{subsec:errorestimator} to obtain well-informed initial step sizes. 

\subsection{Scaled Newton polynomials}\label{subsec:scalednewton}
A by-product of \Cref{alg:adaptsgmres} is that the adaptive block size can serve as a metric for evaluating the effects of other measures on the stability limit of adaptive s-step GMRES. Previous discussions on polynomial bases in literature focus on the conditioning of the Krylov basis matrix $V(Q_i, s)$~\cite{joubert1992parallelizable,joubert1992parallelizable2,philippe2012generation,ballard2014communication,carson2015communication}. However, as we pointed out earlier in \Cref{subsec:stability_QR}, the inter-orthogonalization process could alter the condition of the Krylov matrix, making them ill-conditioned for intra-orthogonalization. Hence, it is insufficient to only monitor the condition number growth of the $V(Q_i, s)$ at the end of MPK. Using the adaptive $s$-step GMRES algorithm, we would like to use the adapted step size $s$ as an additional metric to evaluate different polynomial bases from a stability point of view. 

In addition to the monomial basis introduced in \Cref{eq:monomial}, we are also interested in Newton polynomials. The standard Newton polynomials are defined as
\begin{equation}\label{eq:newton}
    p_j(A) = \Pi_{i=1}^{j}{(A-\theta_i I)}.
\end{equation}
The shifts $\theta_i$ are Ritz values of $A$ in a Leja ordering. In practice, the Ritz values are obtained from $s$ iterations of classical Arnoldi or GMRES first. This is a common strategy for $s$-step GMRES implementations; see, e.g., \cite{bai1994newton,hoemmen2010communication,ballard2014communication} for further details. 

In this work, we introduce an additional variant of the Newton polynomials: \textit{scaled} Newton polynomials defined as
\begin{equation}\label{eq:scalednewton}
    p_j(A) = \Pi_{i=1}^{j}{(A-\theta_i I)/\gamma_i},
\end{equation}
where $\gamma_i$ are the scaling coefficients. The same Leja ordering of Ritz values applies. The standard Newton polynomial is considered a special case with $\gamma_i = 1$. 

In what follows, we propose a heuristic to choose $\gamma_i$. Our analysis will also help to understand the derivation of the initial step size estimator in \Cref{subsec:errorestimator} and the effectiveness of the matrix equilibration in \Cref{sec:test_matrixscaling}.

Let $V_j$ denote the $j$th column of $V=V(Q_i,s)\in\mathbb{R}^{N\times(s+1)}$, which corresponds to the $j$th Krylov basis vector to be orthogonalized. The scaled Newton polynomial can then be written as a two-term recurrence
\begin{equation}\label{eq:newtonrecurrence}
\begin{split}
    V_1 &= Q_i, \\
    V_j &= (AV_{j-1} - \theta_{j-1}V_{j-1})/\gamma_{j-1}.
\end{split}
\end{equation}
Recall that the first vector $V_1 = Q_i$ is an orthonormalized vector at the beginning of each block iteration, and $\theta_i$ are Ritz values that are approximations of true eigenvalues. 

For our analysis, assume that $A$ is diagonalizable
\begin{equation}
    A = U\Lambda U^{-1},
\end{equation}
where $\Lambda$ and $U$ are eigenvalues and eigenvectors of $A$. For simplicity, we assume that all eigenvalues are real, though extensions to eigenvalues with complex conjugate pairs are straightforward. We express $V$ or $V_j$ using eigenvectors of $A$ as
\begin{equation}\label{eq:vectordecomp}
    V = UC, \qquad V_j = UC_j.
\end{equation}
Substituting \Cref{eq:vectordecomp} into \Cref{eq:newtonrecurrence}, entries in the matrix $C$ can be expressed using the recurrence
\[
    C_{i,j+1} = C_{i,j} \frac{\lambda_i - \theta_j}{\gamma_j}.
\]
Given the base case, $V_1 = UC_1$, another way to write the above relationship is
\begin{equation}\label{eq:Cmatrix}
    C_{i,j} = C_{i,1}\prod_{k=1}^{j-1}{\frac{\lambda_i-\theta_k}{\gamma_k}},\qquad 1 \leq i \leq N \text{ and } 1 \leq j \leq s+1.
\end{equation}
This relation describes the evolution of $V_j$ in the eigenspace of $A$ under the action of scaled Newton polynomials. One way to choose the scaling coefficients is to minimize the growth/decay of norm of each Krylov basis vector
\begin{equation}
    \Bigg\lvert\prod_{k=1}^{j-1}{\frac{\lambda_i - \theta_k}{\gamma_k}}\Bigg\rvert \sim O(1).
\end{equation}
Since in practice we do not have access to the true eigenvalue spectrum $\lambda_i$ but only Ritz values $\theta_i$, the scaling coefficients are therefore chosen as
\begin{equation}\label{eq:scalingcoeff}
    \gamma_i = |\bar{\theta} - \theta_i|,
\end{equation}
where $\bar{\theta}$ is the average of all Ritz values. 

Similar to the unscaled Newton polynomials, the Ritz values are obtained from $s$ iterations of classical Arnoldi or GMRES. The average is therefore taken over the $s$ Ritz values computed.

\Cref{eq:scalednewton} and \Cref{eq:scalingcoeff} complete the definition of scaled Newton polynomials. 

\subsection{Initial step size estimator for scaled Newton polynomials}\label{subsec:errorestimator}

As we will demonstrate in numerical experiments in \Cref{sec:tests_stability}, the scaled Newton polynomials can give very large adapted step sizes in many cases, which is desirable for reducing communication latency. However, in certain cases (such as the one in \Cref{sec:test_E20R5000_ILU}), the adapted step size from scaled Newton polynomials can be limited. By blindly setting a very large initial step size for arbitrary matrices, one might end up with relatively small step sizes, resulting in the removal of many columns of the Krylov basis matrix at the first block iteration. The wasted computation in such scenarios often can offset the gain in communication savings. In this section, we show that, using the eigendecomposition analysis in \Cref{subsec:scalednewton}, an estimator can be derived that informs the choice of initial step sizes for scaled Newton polynomials.

\begin{figure}[htbp]
    \centering
    \includegraphics[width=\textwidth]{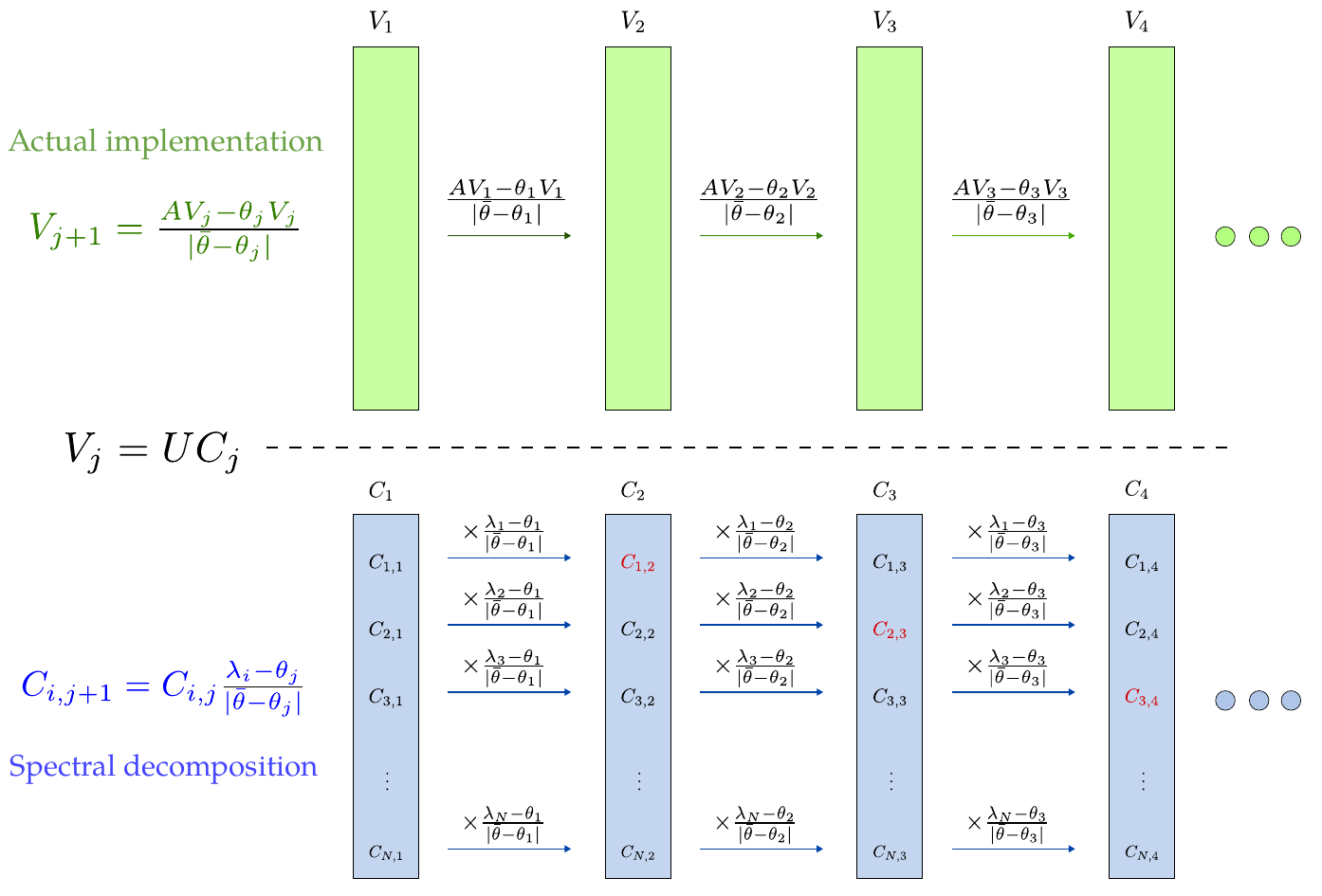}
    \caption{Evolution of Krylov vectors $V_j$ (top) and their corresponding coefficients $C_{i,j}$ in the eigenspace of $A$(bottom) under the action of Scaled Newton polynomials. The two representations are related via $V_j = UC_j$. The entries highlighted in red are $O(\varepsilon)$ if Ritz values are well-approximated.}
    \label{fig:estimator}
\end{figure}

From \Cref{eq:Cmatrix} and \Cref{eq:scalingcoeff}, the generation of Krylov basis vectors using scaled Newton polynomials can be expressed in the eigenspace of $A$ as
\begin{equation}\label{eq:errorindicator}
    C_{i,j} = C_{i,1}\prod_{k=1}^{j-1}{\frac{\lambda_i-\theta_k}{|\bar{\theta}-\theta_k|}},\qquad 1 \leq i \leq N, \; 1 \leq j \leq s+1.
\end{equation}
We are interested in scaled Newton polynomials in this particular form because by using Ritz values to approximate true eigenvalues, we can compute the product term ($\prod_{k=1}^{j-1}{\frac{\lambda_i-\theta_k}{|\bar{\theta}-\theta_k|}}$) approximately before the actual adaptive $s$-step GMRES algorithm. Hence, once Ritz values are computed, this estimator can be computed without any communication cost. For illustration, the evolution of Krylov vectors $V_j$ and their corresponding eigenspace coefficients $C_{i,j}$ are visualized in \Cref{fig:estimator}. Ideally, if the Ritz values $\theta_k$ are exactly equal to the true eigenvalues $\lambda_k$ and roundoff errors are absent, the numerators would cancel out with 
\[
    \lambda_i - \theta_k = 0,  \qquad \forall i = k.
\]
The relation \Cref{eq:errorindicator} for all $1 \leq i \leq N$ can then be reduced to
\[
    C_{i,j} = 
    \begin{cases}
    C_{i,1}\prod_{k=1}^{j-1}{\frac{\lambda_i-\theta_k}{|\bar{\theta}-\theta_k|}},\qquad &1 \leq j \leq i, \\
    0,\qquad &i < j \leq s+1, \\
    \end{cases}
\]
The second case above implies that the strictly upper triangular entries of $C$ are all zero. Now, assuming that some eigenvalues are well-approximated with finite precision, $\lambda_i - \theta_k = O(\varepsilon)$, the second case in the above relation becomes
\begin{equation*}
    C_{i,j} = \Bigg(\prod_{k=1}^{i-1}{\frac{\lambda_i-\theta_k}{|\bar{\theta}-\theta_k|}}\Bigg)O(\varepsilon)\Bigg(\prod_{k=i+1}^{j-1}{\frac{\lambda_i-\theta_k}{|\bar{\theta}-\theta_k|}}\Bigg), \qquad i < j \leq s+1.
\end{equation*}

In other words, the entries on the superdiagonal of $C$ (i.e., $C_{j,j+1}$) do not become strictly zero due to numerical error $O(\varepsilon)$. Such numerical error can either grow or decay in the subsequent application of scaled Newton polynomials depending on the evolution of the product term. Overall, the matrix $C$ is fully dense instead of being lower triangular. 

Recall that the overarching goal of introducing different polynomial bases is to minimize the condition number of the Krylov basis matrix $V$. Although it is difficult to achieve this for an arbitrary starting vector, we can instead heuristically try to optimize for a less stringent metric, the vector-wise norm of Krylov basis vectors, $\|V_j\|_2$. In particular, we want $\|V_j\|_2$ to stay at $O(1)$ ($\|V_1\|_2 = 1$ by construction) since the exponential growth of $\|V_j\|_2$ results in exponential increase in the condition number of the Krylov basis matrix $V$.  In the eigenspace, this translates to $\|C_j\|_2 \approx O(\|C_1\|_2)$. Note that we do not have access to $C_j$ explicitly as the eigenvectors of $A$ are prohibitively expensive to compute, but we merely need to approximate the product term in \Cref{eq:errorindicator} at every $j$ to monitor the growth of $\|C_j\|$. Therefore, we make use of the Ritz values from $s$ iterations of Arnoldi/GMRES to approximate the product term by defining an auxiliary matrix $E_{i,j} \in \mathbb{R}^{s\times s}$
\begin{equation}\label{eq:auxiliarymatrix}
    E_{i,j} = \begin{cases}
    \Pi_{k = 1}^{j-1} \frac{|\theta_i - \theta_k|}{|\bar{\theta} - \theta_k|} & j < i,\\
    \Big(\Pi_{k = 1}^{j-1} \frac{|\theta_i - \theta_k|}{|\bar{\theta} - \theta_k|}\Big) O(\varepsilon) & j = i,\\
    \Big(\Pi_{k = 1}^{i-1} \frac{|\theta_i - \theta_k|}{|\bar{\theta} - \theta_k|}\Big) O(\varepsilon) \Big(\Pi_{k = i+1}^{j-1} \frac{|\theta_i - \theta_k|}{|\bar{\theta} - \theta_k|}\Big) & i < j \leq s.\\
    \end{cases}
\end{equation}
In actual implementation, $O(\varepsilon)$ can be chosen as $\epsilon$, the unit round-off. Each row $i$ of $E_{i,j}$ estimates the evolution of the product term at polynomials of degree $j$ for the coefficients of $C_{i,:}$, including the evolution of the numerical error. The column-wise norm $\|E_{j}\|$ provides an estimate of the ratio $\|C_j\| / \|C_1\|$. Once $\|E_{j}\|$ grows significantly beyond $O(1)$, the scaled Newton polynomials will very likely increase the vector-wise norm of the Krylov basis vectors $V_j$. Since $\|V_1\|_2 = 1$, subsequently increasing the norm of $\|V_j\|_2$ will result in growing incremental condition number of the overall Krylov basis matrix $V$, leading to small step size for stable orthogonalization. 

Let $\Omega_\text{est}$ be a user-specified threshold significantly greater than 1. We propose a way to estimate the initial step size via 
\begin{equation}
    s_0^* = \argmax_j \{\|E_{j}\| < \Omega_\text{est}\} 
\end{equation}
where $s_0^*$ denotes the maximum initial step size predicted. As the Ritz values are available on all MPI ranks, the estimator requires no communication and can be used to minimize the number of columns that partial CholQR has to truncate at the \textit{first} block iteration. In our experience, a threshold of $\Omega_\text{est} = 10^{-1}\epsilon^{-1/2}$ in for the estimator does a good job at limiting the initial step size and gives reasonable prediction $s_0^*$, though more conservative values can be used.

\section{Results on numerical stability}\label{sec:tests_stability}

We numerically verify the proposed adaptive $s$-step GMRES algorithm to demonstrate the convergence and stability of the algorithm. In this section, we implement \Cref{alg:adaptsgmres} in MATLAB and use double precision for all tests where $\epsilon \approx 2^{-53}$. As the work here focuses on the algorithm's stability, it suffices to verify the implementation sequentially. We assume the corresponding parallel implementations only introduce additional data movements and do not affect the stability of the overall algorithm. Benchmarking the algorithmic scalability in a distributed network will be presented in \Cref{sec:scalability}.

The metrics we are interested in are the 2-norm of the relative residual at iteration $i$,
\begin{equation}
    \|b - Ax_i\|_2/\|b-Ax_0\|_2,
\end{equation}
as well as the LOO error of the block QR scheme,
\begin{equation}
    \|I - Q_{1:i}^TQ_{1:i}\|_F,
\end{equation}
where $Q$ contains the orthogonal basis vectors. 

We are also interested in the adapted block size $s$ since it can serve as a metric for evaluating the effects of other measures on the stability limit of adaptive $s$-step GMRES. 

In all tests, an initial step size $s_0$ is set at the beginning of the numerical experiment. The stopping criteria for partial CholQR are imposed at $\Omega = 10^7$. ICE is used for condition number bounds and is verified against SVD estimates (i.e., \texttt{cond} in MATLAB).

\subsection{Numerical stability of the base algorithm}
We start with the simplest form of the algorithm in this subsection. We use a \textit{monomial} basis defined by \Cref{eq:monomial} in the MPK, and no preconditioners are used to condition the linear systems. Although convergence tends to be slow in the absence of preconditioners for general problems, we place our emphasis on the agreement between the standard GMRES algorithm and the adaptive GMRES algorithm. Numerical tests with different polynomial bases and preconditioning techniques are delayed to \Cref{sec:tests_polynomialbasis} and \Cref{sec:tests_preconditioners}.

\subsubsection{Diagonal matrix}\label{subsubsec:test_diagonalmatrix}

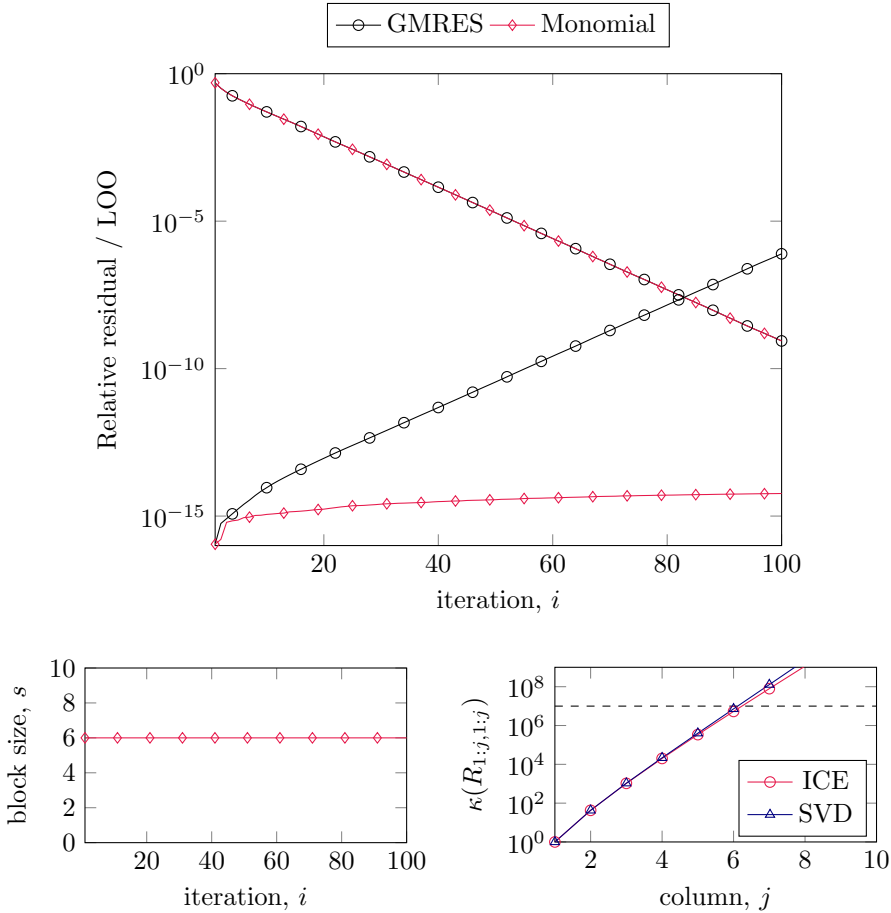
\begin{figure}[tbhp]
\centering
\subfloat{
    \begin{tikzpicture}
    \begin{axis}[
    width=0.7\textwidth,
    ymode=log,
    xmin=1,
    xmax=100,
    ytick={1e-15, 1e-10, 1e-5, 1},
    ymax=1,
    ymin=1e-16,
    xlabel={iteration, $i$}, 
    ylabel={Relative residual / LOO},
    legend entries={GMRES, Monomial},
    % legend pos=outer north east,
    legend columns=-1,
    legend style={at={(0.5,1.15)},anchor=north},
    ]
    \addplot [mark=o,mark repeat=6,mark phase=4] table [x=i,y=res] {Figures_data/Diagonal/Diagonal_res_LOO_s.dat};
    \addplot [
    hcred,
    mark=diamond,
    mark repeat=6,
    ] table [x=i,y=res_CA] {Figures_data/Diagonal/Diagonal_res_LOO_s.dat};
    \addplot [mark=o,mark repeat=6,mark phase=4] table [x=i,y=LOO] {Figures_data/Diagonal/Diagonal_res_LOO_s.dat};
    \addplot [
    hcred, 
    mark=diamond,
    mark repeat=6,
    ] table [x=i,y=LOO_CA] {Figures_data/Diagonal/Diagonal_res_LOO_s.dat};
    \end{axis}
    \end{tikzpicture}
}
\\
\subfloat{
    \begin{tikzpicture}
    \begin{axis}[
    width=0.45\textwidth,
    height=0.3\textwidth,
    xmin=1,
    xmax=100,
    ymin=0,
    ymax=10,
    xlabel={iteration, $i$}, 
    ylabel={block size, $s$},
    ]
    \addplot [
    hcred,
    mark=diamond,
    mark repeat=10,
    ] table [x=i,y=s] {Figures_data/Diagonal/Diagonal_res_LOO_s.dat};
    \end{axis}
    \end{tikzpicture}
}
\subfloat{
    \begin{tikzpicture}
    \begin{semilogyaxis}[
    width=0.45\textwidth,
    height=0.3\textwidth,
    xmin=1,
    xmax=10,
    ymin=1,
    ymax=1e9,
    xlabel={column, $j$}, 
    ylabel={$\kappa(R_{1:j,1:j})$},
    ytick={1e0, 1e2, 1e4, 1e6, 1e8},
    legend entries={ICE, SVD},
    legend pos=south east,
    ]
    \addplot [
    hcred, 
    mark=o,
    ] table [x=i,y=ICE] {Figures_data/Diagonal/Diagonal_cond.dat};
    \addplot [
    hcnavy,
    mark=triangle,
    ] table [x=i,y=SVD] {Figures_data/Diagonal/Diagonal_cond.dat};
    \addplot [
    sharp plot, dashed,
    ] coordinates
    {(1,1e7) (10,1e7)};
    \end{semilogyaxis}
    \end{tikzpicture}
}
\caption{Diagonal matrix of size $N=10^4$ with evenly distributed eigenvalues in (0.1,10). The top sub-figure shows the relative residual and the LOO errors of GMRES and the adaptive $s$-step GMRES with a monomial basis. The relative residual of the adaptive $s$-step GMRES algorithm agrees with the baseline. The LOO error of the adaptive s-step GMRES algorithm is maintained near the machine epsilon. The lower left sub-figure shows the adapted block size, $s$. In this case, the adapted step size is constant at $s = 6$. The lower right sub-figure shows the incremental condition number of the triangular matrix seen by the first instance of partial CholQR in the first block iteration. The ICE estimate is very close to the true value given by the SVD.}
\label{fig:diagonalmatrix}
\end{figure}

To demonstrate the stability and convergence of the algorithm, we first use a diagonal matrix of size $N=10^4$ with evenly distributed eigenvalues in $(0.1, 10)$. The initial step size is $s_0 = 10$. The top sub-figure in \Cref{fig:diagonalmatrix} shows the relative residual and LOO error for both the standard GMRES and the adaptive $s$-step GMRES algorithm (with a monomial basis). Here, the relative residual curves start close to 1 and gradually decrease as the solutions converge, whereas the LOO curves start close to machine epsilon and increase at different rates depending on the orthogonalization scheme in the algorithm. The convergence of the adaptive algorithm is in good agreement with the GMRES baseline. The LOO of GMRES increases gradually due to the Modified Gram-Schmidt (MGS) orthogonalization used in the standard Arnoldi procedure. The LOO error of the proposed algorithm stays at $O(\epsilon)$. Note that the LOO error associated with BCGS2 with CholQR for the same matrix and monomial basis was also analyzed in Fig.~5 of~\cite{carson2022block}. Their numerical experiments show a significantly higher LOO error because the authors chose a constant step size of $s=10$, which is above the stability limit of the QR scheme. The lower left sub-figure in \Cref{fig:diagonalmatrix} shows the block size at each iteration. The larger the block size, the fewer block iterations are used for the same number of iterations. Therefore, the ``area under the graph'' for the block size plot indicates the amount of potential communication savings. With an initial user-defined step size of $s_0=10$, the algorithm adapted the step size to $s=6$. If we plot the incremental condition number growth of the triangular matrix (i.e., $R_{1:j,1:j}$ in \Cref{alg:pcholICE}) for the first partial CholQR instance in the first block iteration, we also observe that $\kappa(R_{1:j,1:j}) > O(\epsilon^{-1/2})$ after $s=6$, indicating that larger value of $s$ would become unstable for intra-orthogonalization. Therefore, the algorithm has to truncate the generated Krylov basis matrix down to $s=6$ at the first block iteration and continue with this step size to maintain stability in the orthogonalization process. We also point out that in \Cref{fig:diagonalmatrix}, ICE gives reasonably accurate estimates of the incremental condition number compared to the more expensive SVD decomposition. In general, we find that $\Omega = 10^7$ is usually sufficient to bound the condition number in double precision arithmetic.

 %% lengthy figures are isolated into individual files

\subsubsection{2D Laplace}\label{sec:test_2Dlaplace}

\begin{figure}[tbhp]
\centering
\subfloat{
    \begin{tikzpicture}
    \begin{groupplot}[
        group style={
            group name=my fancy plots,
            group size=1 by 2,
            xticklabels at=edge bottom,
            vertical sep=0pt,
        },
        width=0.7\textwidth,
        xmin=0, xmax=500,
    ]
    \nextgroupplot[
        ymode=log,
        ymin=0.05,ymax=1,
        ytick={1e-1,1},
        axis x line=top, 
        axis y discontinuity=parallel,
        height=0.35\textwidth,legend entries={GMRES, Monomial},
        % legend pos=outer north east,
        legend columns=-1,
        legend style={at={(0.5,1.3)},anchor=north},
    ]
    \addplot [mark=o,mark repeat=26,mark phase=14] table [x=i,y=res] {Figures_data/2DLaplace/2DLaplace_res_LOO_s.dat};     
    \addplot [
    hcred,
    mark=diamond,
    mark repeat=26,
    ] table [x=i,y=res_CA] {Figures_data/2DLaplace/2DLaplace_res_LOO_s.dat};   
    
    \nextgroupplot[
        ymode=log,
        ymin=1e-16,ymax=1e-10,
        ytick={1e-15,1e-13,1e-11},
        axis x line=bottom,
        x axis line style={},
        % x axis line style={-stealth},
        height=0.35\textwidth,
        xlabel={iteration, $i$},
        ylabel={Relative residual / LOO},
        ylabel style={at={(ticklabel cs:1)}},
    ]
    \addplot [mark=o,mark repeat=26,mark phase=14] table [x=i,y=LOO] {Figures_data/2DLaplace/2DLaplace_res_LOO_s.dat};
    \addplot [
    hcred, 
    mark=diamond,
    mark repeat=26,
    ] table [x=i,y=LOO_CA] {Figures_data/2DLaplace/2DLaplace_res_LOO_s.dat};
    \end{groupplot}
    \end{tikzpicture}
}
\\
\subfloat{
    \begin{tikzpicture}
    \begin{axis}[
    width=0.45\textwidth,
    height=0.3\textwidth,
    xmin=1,
    xmax=500,
    ymin=0,
    ymax=10,
    xlabel={iteration, $i$}, 
    ylabel={block size, $s$},
    xtick={250,500},
    ]
    \addplot [
    hcred, 
    mark=diamond,
    mark repeat=50,
    ] table [x=i,y=s] {Figures_data/2DLaplace/2DLaplace_res_LOO_s.dat};
    \end{axis}
    \end{tikzpicture}
}
\subfloat{
    \begin{tikzpicture}
    \begin{semilogyaxis}[
    width=0.45\textwidth,
    height=0.3\textwidth,
    xmin=1,
    xmax=10,
    ymin=1,
    ymax=1e9,
    xlabel={column, $j$}, 
    ylabel={$\kappa(R_{1:j,1:j})$},
    ytick={1e0, 1e2, 1e4, 1e6, 1e8},
    legend entries={ICE, SVD},
    legend pos=south east,
    ]
    \addplot [
    hcred, 
    mark=o,
    ] table [x=i,y=ICE] {Figures_data/2DLaplace/2DLaplace_cond.dat};
    \addplot [
    hcnavy,
    mark=triangle,
    ] table [x=i,y=SVD] {Figures_data/2DLaplace/2DLaplace_cond.dat};
    \addplot [
    sharp plot, dashed,
    ] coordinates
    {(1,1e7) (10,1e7)};
    \end{semilogyaxis}
    \end{tikzpicture}
}
\caption{2D Laplace matrix of size $N=1.6\times10^\rpy{5}$ using 5-stencil finite difference discretization. The top sub-figure shows the relative residual and LOO errors of GMRES and adaptive s-step GMRES with a monomial basis. Good agreement of the relative residuals between the baseline and the adaptive algorithm can be observed. LOO of the adaptive s-step GMRES is kept near the machine epsilon. The lower left sub-figure shows the adapted block size, $s$. In this case, the adapted step size is constant at $s=6$. The lower right sub-figure shows the incremental condition number of the triangular matrix seen by the first instance of partial CholQR in the first block iteration. The ICE estimate agrees well with the true SVD estimate.}
\label{fig:2Dlaplace}
\end{figure}
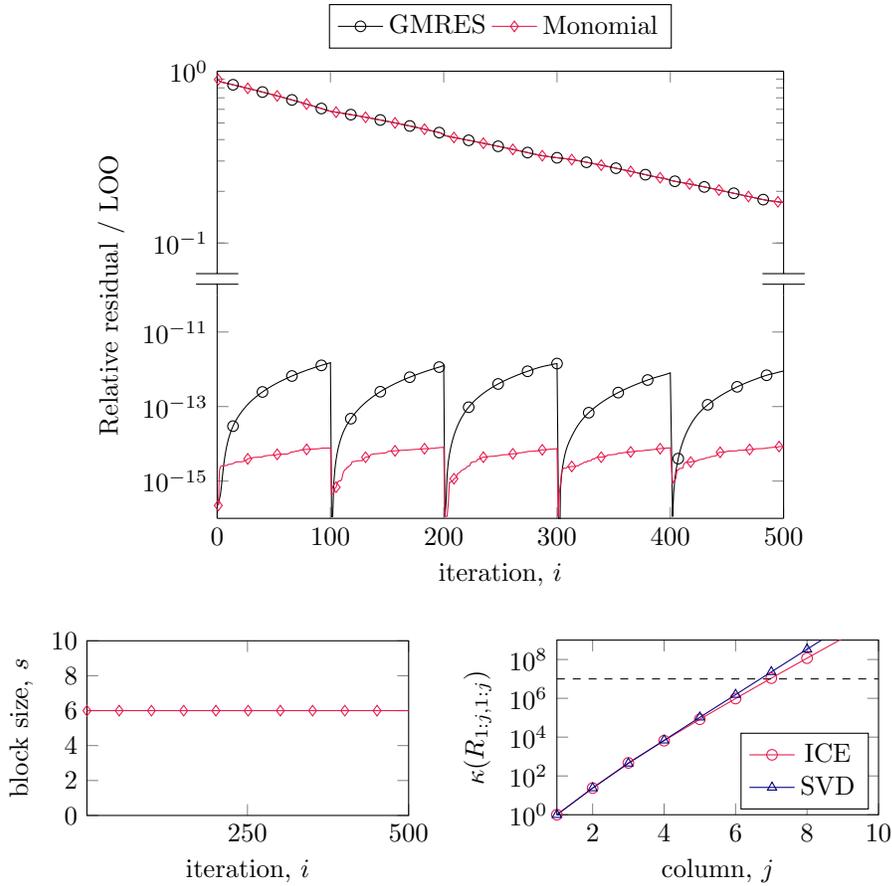

The second matrix comes from solving the standard 2D Laplace problem on a $400^2$ uniform grid using a 5-stencil finite difference discretization, resulting in a matrix dimension of $N=1.6\times10^5$. A maximum inner iteration of 100 with 5 restarts is used. Without any preconditioner, the relative residual slowly converges. However, excellent agreement for residuals is still observed in \Cref{fig:2Dlaplace}. LOO of adaptive $s$-step GMRES is kept at $O(\epsilon)$. The periodic pattern of LOO is a consequence of the reconstruction of orthogonal bases at the beginning of each restart. 

It should be noted that since the maximum inner iteration is capped at 100 (which is not divisible by $s=6$ as given by \Cref{fig:2Dlaplace}), some Krylov basis vectors in the last block are discarded, even though they do not undermine the stability of the algorithm. This is done only for benchmarking purposes, since we are comparing this with the GMRES baseline. In practice, one is free to extend the maximum inner iterations within each restart to retain the entire last block of Krylov basis vectors for better convergence and communication savings.

  %% lengthy figures are isolated into individual files

\subsection{Numerical experiments with different polynomial bases}\label{sec:tests_polynomialbasis}
In this subsection, we compare the monomial basis \Cref{eq:monomial}, the standard Newton polynomial basis \Cref{eq:newton} as well as the scaled Newton polynomial basis \Cref{eq:scalednewton} with the baseline GMRES. The primary emphasis is placed on the adapted step size for various bases.

\subsubsection{Diagonal matrix}
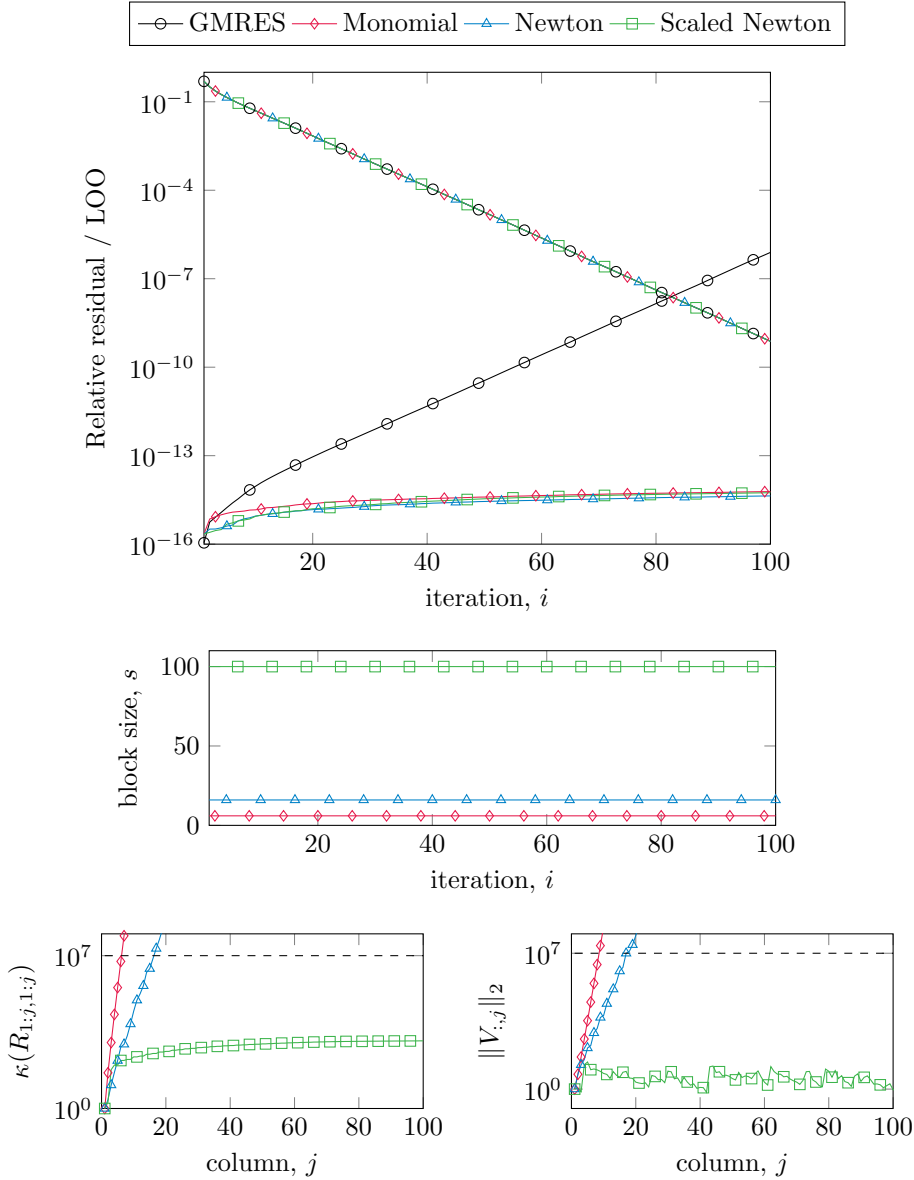
\begin{figure}[!tbhp]
\centering
\subfloat{
    \begin{tikzpicture}
    \begin{axis}[
    width=0.7\textwidth,
    ymode=log,
    xmin=1,
    xmax=100,
    ymin=1e-16,
    ymax=1,
    xlabel={iteration, $i$}, 
    ylabel={Relative residual / LOO},
    legend entries={GMRES, Monomial, Newton, Scaled Newton},
    legend columns=-1,
    legend style={at={(0.5,1.15)},anchor=north},
    ]
    \addplot [mark=o,mark repeat=8] table [x=i,y=GMRES] {Figures_data/Diagonal_Poly/Polynomial_res.dat};
    \addplot [
    hcred, 
    mark=diamond,
    mark repeat=8,
    mark phase=3,
    ] table [x=i,y=Monomial] {Figures_data/Diagonal_Poly/Polynomial_res.dat};
    \addplot [
    hcblue, 
    mark=triangle,
    mark repeat=8,
    mark phase=5,
    ] table [x=i,y=Newton] {Figures_data/Diagonal_Poly/Polynomial_res.dat};
    \addplot [
    hcgreen, 
    mark=square,
    mark repeat=8,
    mark phase=7,
    ] table [x=i,y=SNewton] {Figures_data/Diagonal_Poly/Polynomial_res.dat};
    % \addplot [
    % hcgrey,
    % mark=o,
    % mark repeat=8,
    % mark phase=8,
    % ] table [x=i,y=Chebyshev] {Figures_data/Diagonal_Poly/Polynomial_res.dat};
    \addplot [mark=o,mark repeat=8] table [x=i,y=LOO] {Figures_data/Diagonal/Diagonal_res_LOO_s.dat};
    \addplot [
    hcred, 
    mark=diamond,
    mark repeat=8,
    mark phase=3,
    ] table [x=i,y=Monomial] {Figures_data/Diagonal_Poly/Polynomial_LOO.dat};
    \addplot [
    hcblue, 
    mark=triangle,
    mark repeat=8,
    mark phase=5,
    ] table [x=i,y=Newton] {Figures_data/Diagonal_Poly/Polynomial_LOO.dat};
    \addplot [
    hcgreen, 
    mark=square,
    mark repeat=8,
    mark phase=7,
    ] table [x=i,y=SNewton] {Figures_data/Diagonal_Poly/Polynomial_LOO.dat};
    % \addplot [
    % hcgrey, 
    % mark=o,
    % mark repeat=8,
    % mark phase=8,
    % ] table [x=i,y=Chebyshev] {Figures_data/Diagonal_Poly/Polynomial_LOO.dat};
    \end{axis}
    \end{tikzpicture}
}
\\
\subfloat{
    \begin{tikzpicture}
    \begin{axis}[
    width=0.7\textwidth,
    height=0.3\textwidth,
    xmin=1,
    xmax=100,
    ymin=0,
    ymax=110,
    xlabel={iteration, $i$}, 
    ylabel={block size, $s$},
    ]
    \addplot [
    hcred, 
    mark=diamond,
    mark repeat=6,
    mark phase=2,
    ] table [x=i,y=Monomial] {Figures_data/Diagonal_Poly/Polynomial_s.dat};
    \addplot [
    hcblue, 
    mark=triangle,
    mark repeat=6,
    mark phase=4,
    ] table [x=i,y=Newton] {Figures_data/Diagonal_Poly/Polynomial_s.dat};
    \addplot [
    hcgreen, 
    mark=square,
    mark repeat=6,
    mark phase=6,
    ] table [x=i,y=SNewton] {Figures_data/Diagonal_Poly/Polynomial_s.dat};
    % \addplot [
    % hcgrey, 
    % mark=o,
    % mark repeat=8,
    % mark phase=8,
    % ] table [x=i,y=Chebyshev] {Figures_data/Diagonal_Poly/Polynomial_s.dat};
    \end{axis}
    \end{tikzpicture}
}\;
\\
\subfloat{
    \begin{tikzpicture}
    \begin{axis}[
    ymode=log,
    width=0.45\textwidth,
    height=0.3\textwidth,
    xmin=0,
    xmax=100,
    ymin=1,
    ymax=1e8,
    ytick={1,1e7},
    xlabel={column, $j$}, 
    ylabel={$\kappa(R_{1:j,1:j})$},
    ]
    \addplot [
    hcred,
    mark=diamond,
    ] table [x=i,y=Monomial_ICE] {Figures_data/Diagonal_Poly/Polynomial_cond.dat};
    \addplot [
    hcblue,
    mark=triangle,
    mark repeat=2,
    ] table [x=i,y=Newton_ICE] {Figures_data/Diagonal_Poly/Polynomial_cond.dat};
    \addplot [
    hcgreen,
    mark=square,
    mark repeat=5,
    ] table [x=i,y=SNewton_ICE] {Figures_data/Diagonal_Poly/Polynomial_cond.dat};
    % \addplot [
    % hcgrey,
    % mark=o,
    % mark repeat=2,
    % ] table [x=i,y=Chebyshev_ICE] {Figures_data/Diagonal_Poly/Polynomial_cond.dat};
    \addplot [
    sharp plot, dashed,
    ] coordinates
    {(1,1e7) (100,1e7)};
    \end{axis}
    \end{tikzpicture}
}
\subfloat{
    \begin{tikzpicture}
    \begin{axis}[
    ymode=log,
    width=0.45\textwidth,
    height=0.3\textwidth,
    xmin=0,
    xmax=100,
    ymin=1e-1,
    ymax=1e8,
    ytick={1,1e7},
    xlabel={column, $j$}, 
    ylabel={$\|V_{:,j}\|_2$},
    ]
    \addplot [
    hcred,
    mark=diamond,
    ] table [x=i,y=monomial_v] {Figures_data/Diagonal_Poly/Polynomial_vecnorm.dat};
    \addplot [
    hcblue,
    mark=triangle,
    mark repeat=2,
    ] table [x=i,y=newton_v] {Figures_data/Diagonal_Poly/Polynomial_vecnorm.dat};
    \addplot [
    hcgreen,
    mark=square,
    mark repeat=5,
    ] table [x=i,y=snewton_v] {Figures_data/Diagonal_Poly/Polynomial_vecnorm.dat};
    % \addplot [
    % hcgrey,
    % mark=o,
    % mark repeat=2,
    % ] table [x=i,y=chebyshev_v] {Figures_data/Diagonal_Poly/Polynomial_vecnorm.dat};
    \addplot [
    sharp plot, dashed,
    ] coordinates
    {(1,1e7) (100,1e7)};
    \end{axis}
    \end{tikzpicture}
}
\caption{Diagonal matrix. Top figure: relative residual and LOO errors of GMRES and adaptive $s$-step GMRES using different polynomial bases in the MPK. Middle figure: adapted step size $s$. The scaled Newton polynomials have the largest adapted step size of $s=100$. Lower left figure: incremental condition number of the triangular matrix seen by the first instance of partial CholQR in the first block iteration. Lower right figure: vector-wise norm of the Krylov basis matrix. The scaled Newton polynomial is effective at keeping the vector norm to $O(1)$ and slowing the incremental condition number growth.}
\label{fig:diagonalmatrix_poly}
\end{figure}

We use the same diagonal matrix as \Cref{subsubsec:test_diagonalmatrix} with different polynomial bases. The initial step size is $s_0 = 100$ to illustrate the stability limit associated with each polynomial basis. The results are summarized in \Cref{fig:diagonalmatrix_poly}. All polynomial bases give stable results that agree well with the baseline GMRES. In terms of step sizes, all polynomial bases outperform the monomial basis. We would like to highlight that the scaled Newton polynomial stands out because the adapted step size is the same as the number of iterations $s = 100$. In other words, the implementation using scaled Newton polynomials requires only \textit{one} block iteration for the problem to converge. This is because of the much better conditioned Krylov basis matrix, as shown in the lower left sub-figure of \Cref{fig:diagonalmatrix_poly}. The difference in the adapted step sizes is reflected in the vector-wise norm of the Krylov basis matrix. We observe that the scaled Newton polynomials are capable of keeping the column vector norms to $O(1)$ for a large number of iterations, thus helping to slow down the condition number growth of the Krylov basis matrix.

  %% lengthy figures are isolated into individual files

\subsubsection{E20R5000 matrix}\label{sec:test_E20R5000}
The next test matrix is E20R5000 from Matrix Market~\cite{boisvert1997matrix}, which models 2D fluid flow in a driven cavity. The matrix is known to be difficult for Krylov iterative solvers. In the absence of any preconditioning technique, the relative residual decreases very slowly, as shown in \Cref{fig:E20R5000}. However, the various implementations agree with the baseline GMRES algorithm, and the LOO errors are near machine accuracy. With an initial step size of $s_0=150$, the scaled Newton polynomial allows the largest stable step size among all polynomial bases. However, there is a sharp decrease in the step size after the first block iteration. To understand this, we further plot the incremental condition number of both the input block to partial CholQR ($\kappa(R_{1:j,1:j})$) and the Krylov basis matrix ($\kappa(V_{:,1:j})$) in the two lower subfigures of \Cref{fig:E20R5000}. We observe that the Krylov basis matrix is reasonably well-conditioned in both the first and second iterations, but the input block to partial CholQR becomes ill-conditioned rapidly in the second block iteration. This is attributed to the inter-orthogonalization process. The Krylov basis matrix in the second block iteration has to be inter-orthogonalized with respect to the first block of size $s\approx 100$ before partial CholQR factorization. The update in the inter-orthogonalization step altered the Krylov basis matrix and caused the condition number to grow exponentially. This demonstrates that a good choice of polynomial basis in MPK does not necessarily guarantee that the orthogonalization is stable for a fixed step size. 

In practice, especially on large distributed machines, one could use a smaller restart length and more restarts to trade convergence for communication savings. If one is satisfied by the large adapted step size of the first block iteration, restarting immediately after one block iteration gives a minimal number of global synchronizations for a fixed number of total iterations. For instance, as shown in \Cref{fig:E20R5000}, adjusting the restart length to $s\approx 100$ reduces the number of block iterations to one in each restart. Although this increases the number of restart cycles from 4 to 6, the total number of global reductions is reduced, leading to communication savings.

\begin{figure}[!tbhp]
\centering
\subfloat{
    \begin{tikzpicture}
    \begin{groupplot}[
        group style={
            group name=my fancy plots,
            group size=1 by 2,
            xticklabels at=edge bottom,
            vertical sep=0pt,
        },
        width=0.7\textwidth,
        xmin=0, xmax=600,
    ]
    \nextgroupplot[
        ymode=log,
        ymin=0.92,ymax=1,
        ytick={0.95,1},
        axis x line=top, 
        axis y discontinuity=parallel,
        height=0.3\textwidth,
        legend entries={GMRES,Monomial,Newton,Scaled Newton},
        legend columns=-1,
        legend style={at={(0.42,1.35)},anchor=north},
    ]
    \addplot [mark=o,mark repeat=40] table [x=i,y=GMRES] {Figures_data/E20R5000/E20R5000_res.dat};
    \addplot [
    hcred,
    mark=diamond,
    mark repeat=40,
    mark phase=11,
    ] table [x=i,y=Monomial] {Figures_data/E20R5000/E20R5000_res.dat};   
    \addplot [
    hcblue,
    mark=triangle,
    mark repeat=40,
    mark phase=21,
    ] table [x=i,y=Newton] {Figures_data/E20R5000/E20R5000_res.dat};
    \addplot [
    hcgreen,
    mark=square,
    mark repeat=40,
    mark phase=31,
    ] table [x=i,y=SNewton] {Figures_data/E20R5000/E20R5000_res.dat};
    % \addplot [
    % hcgrey,
    % mark=o,
    % mark repeat=40,
    % mark phase=40,
    % ] table [x=i,y=Chebyshev] {Figures_data/E20R5000/E20R5000_res.dat};
    
    \nextgroupplot[
        ymode=log,
        ymin=1e-16,ymax=1e-13,
        ytick={1e-16,1e-14},
        xtick={150,300,450,600},
        axis x line=bottom,
        height=0.3\textwidth,
        xlabel={iteration, $i$},
        ylabel={Relative residual / LOO},
        ylabel style={at={(ticklabel cs:1)}},
    ]
    \addplot [mark=o,mark repeat=40] table [x=i,y=GMRES] {Figures_data/E20R5000/E20R5000_LOO.dat};
    \addplot [
    hcred, 
    mark=diamond,
    mark repeat=30,
    mark phase=11,
    ] table [x=i,y=Monomial] {Figures_data/E20R5000/E20R5000_LOO.dat};     
    \addplot [
    hcblue, 
    mark=triangle,
    mark repeat=30,
    mark phase=21,
    ] table [x=i,y=Newton] {Figures_data/E20R5000/E20R5000_LOO.dat};  
    \addplot [
    hcgreen, 
    mark=square,
    mark repeat=30,
    mark phase=31,
    ] table [x=i,y=SNewton] {Figures_data/E20R5000/E20R5000_LOO.dat};  
    % \addplot [
    % hcgrey, 
    % mark=o,
    % mark repeat=40,
    % mark phase=40,
    % ] table [x=i,y=Chebyshev] {Figures_data/E20R5000/E20R5000_LOO.dat};  
    \end{groupplot}
    \end{tikzpicture}
}
\\\,\,
\subfloat{
    \begin{tikzpicture}
    \begin{axis}[
    width=0.7\textwidth,
    height=0.3\textwidth,
    xmin=1,
    xmax=600,
    xtick={150,300,450,600},
    ymin=0,
    ymax=120,
    ytick={0, 40, 80, 120},
    xlabel={iteration, $i$}, 
    ylabel={block size, $s$},
    ylabel shift=5,
    ]
    \addplot [
    hcred,
    mark=diamond,
    mark repeat=30,
    mark phase=15,
    ] table [x=i,y=Monomial] {Figures_data/E20R5000/E20R5000_s.dat};
    \addplot [
    hcblue, 
    mark=triangle,
    mark repeat=30,
    ] table [x=i,y=Newton] {Figures_data/E20R5000/E20R5000_s.dat};
    \addplot [
    hcgreen, 
    mark=square,
    mark repeat=30,
    ] table [x=i,y=SNewton] {Figures_data/E20R5000/E20R5000_s.dat};
    % \addplot [
    % hcgrey, 
    % mark=o,
    % mark repeat=30,
    % ] table [x=i,y=Chebyshev] {Figures_data/E20R5000/E20R5000_s.dat};
    \end{axis}
    \end{tikzpicture}
}
\\
\subfloat{
    \begin{tikzpicture}
    \begin{axis}[
    ymode=log,
    width=0.45\textwidth,
    height=0.4\textwidth,
    xmin=0,
    xmax=150,
    ymin=1,
    ymax=1e8,
    ytick={1,1e7},
    xlabel={column, $j$}, 
    ylabel={$\kappa(R_{1:j,1:j})$},
    legend entries={iter. 1, iter. 2},
    legend pos=south east,
    ]
    \addplot [
    hcgreen, 
    mark=square*,
    mark repeat=8,
    ] table [x=i,y=SVD_restart1] {Figures_data/E20R5000/E20R5000_cond.dat};
    \addplot [
    hcgreen, 
    mark=square,
    mark repeat=2,
    ] table [x=i,y=SVD_restart2] {Figures_data/E20R5000/E20R5000_cond.dat};
    \addplot [
    sharp plot, dashed,
    ] coordinates
    {(1,1e7) (150,1e7)};
    \end{axis}
    \end{tikzpicture}
}
\subfloat{
    \begin{tikzpicture}
    \begin{axis}[
    ymode=log,
    width=0.45\textwidth,
    height=0.4\textwidth,
    xmin=0,
    xmax=150,
    ymin=1,
    ymax=1e8,
    ytick={1,1e7},
    xlabel={column, $j$}, 
    ylabel={$\kappa(V_{:,1:j})$},
    legend entries={iter. 1, iter. 2},
    legend pos=south east,
    ]
    \addplot [
    hcgreen, 
    mark=square*,
    mark repeat=8,
    ] table [x=i,y=SVD1] {Figures_data/E20R5000/E20R5000_condBasis.dat};
    \addplot [
    hcgreen, 
    mark=square,
    mark repeat=8,
    ] table [x=i,y=SVD2] {Figures_data/E20R5000/E20R5000_condBasis.dat};
    \addplot [
    sharp plot, dashed,
    ] coordinates
    {(1,1e7) (150,1e7)};
    \end{axis}
    \end{tikzpicture}
}
\caption{E20R5000 matrix. Top figure: relative residual and LOO errors of GMRES and adaptive $s$-step GMRES using different polynomial bases in the MPK. Middle figure: adapted step size $s$. Scaled Newton polynomials have the largest adapted step size of $s\approx100$. Lower left figure: incremental condition number of the triangular matrix seen by the first instance of partial CholQR for the first two block iterations. Lower right figure: incremental condition number of the Krylov basis matrix for the first two block iterations. Only the scaled Newton polynomial is shown in the last two figures. Despite the well-conditioned Krylov basis matrix in the second iteration, the triangular matrix seen by partial CholQR becomes ill-conditioned rapidly due to the inter-orthogonalization update.}
\label{fig:E20R5000}
\end{figure}  %% lengthy figures are isolated into individual files

\subsection{Numerical experiments with different preconditioning techniques}\label{sec:tests_preconditioners}

Preconditioners can be easily incorporated into \Cref{alg:adaptsgmres}, and the analysis presented so far for generic matrices extends to preconditioned matrices. We first revisit a particular class of preconditioning techniques, matrix equilibration, which is frequently used in the $s$-step GMRES literature for stability considerations rather than to speed up convergence. We numerically investigate the effectiveness of matrix equilibration on conditioning Krylov basis matrices. Subsequently, we present the numerical stability of the adaptive $s$-step GMRES for a more general preconditioning technique, incomplete LU factorization (ILU).

\subsubsection*{A short review on matrix equilibration}\label{sec:test_matrixscaling}

Matrix equilibration has been commonly used in the $s$-step GMRES literature~\cite{hoemmen2010communication,ballard2014communication, carson2015communication} as a preconditioning technique to normalize the matrix. By scaling the spectral radius to $O(1)$, matrix equilibration has been found to be effective in conditioning Krylov basis matrices for different polynomial bases. Ballard et al.\ relied entirely on matrix equilibration techniques and found that standard Newton polynomials without any scaling are sufficient (Section 8.5.1 of~\cite{ballard2014communication}). With the introduction of \textit{scaled} Newton polynomials in this work, we reinvestigate matrix equilibration methods for Newton polynomials.

Matrix equilibration aims to scale the rows and/or columns of the matrix in order to prevent the rapid growth/decay of condition number of Krylov basis vectors. In general, matrix equilibration replaces matrix $A$ by
\begin{equation}\label{eq:equil}
    A' = D_r A D_c,
\end{equation}
with two diagonal matrices $D_r$ and $D_c$. There are many equilibration strategies. Here, we consider two approaches: scalar scaling and column scaling. Scalar scaling scales the entire matrix by a constant (e.g., $D_r = D_c = \sqrt{\alpha}I$) without modifying the relative position of eigenvalues. The constant is often chosen to be related to the spectral radius of $A$ to normalize the spectrum. Column scaling scales each column of the matrix (e.g., $D_r = I$), leading to a modified spectrum of magnitude $O(1)$. Both scaling methods do not require any form of communication. Column scaling can be efficiently computed if the sparse matrix is stored in Compressed Sparse Column (CSC) format. (The same argument applies to row scaling if the sparse matrix is stored in Compressed Sparse Row (CSR) format instead.)

\subsubsection{E20R5000 matrix with matrix equilibration}\label{sec:test_E20R5000_matrixequil}
The E20\-R5000 matrix problem of \Cref{sec:test_E20R5000} is preconditioned with the two matrix equilibration techniques. The spectral radius of $A$ is used for the scalar scaling, i.e., $\alpha = \max{|\theta_i|}$. In both cases, the relative residual plots show good agreement with the baseline GMRES and the LOO is minimal. However, comparing the step size adapted in \Cref{fig:E20R5000_scaling} with those in \Cref{fig:E20R5000}, we observe different levels of improvement for the two Newton polynomials with different scaling methods. The observation can be explained by drawing an analogy to the approximation theory. As noted in~\cite{beckermann2000condition, hoemmen2010communication}, the generation of the Krylov basis using polynomials is equivalent to a polynomial interpolation with the eigenvalues of $A$ being the interpolation nodes. The quality of the interpolation nodes is directly related to the incremental condition number growth of Krylov basis matrices. 

The eigenvalue spectra for the original matrix and preconditioned matrices are plotted in \Cref{fig:E20R5000_spectrum}. Scalar scaling effectively scales the interpolation interval to $O(1)$. Hence, monomial and standard Newton polynomial bases benefit from such normalization, and the adapted step sizes increase. The scaled Newton polynomials see no improvement for scalar scaling because of the scaling coefficient chosen in \Cref{eq:scalingcoeff}. In other words, scaled Newton polynomials are insensitive to scalar scalings. Note that even with a normalized spectrum, there is still a significant difference between the adapted step sizes of two Newton polynomials under scalar scaling. The difference comes from the difference in the scaling coefficient. Standard Newton polynomial under scalar scaling is equivalent to a scalar-scaled Newton polynomial with a constant scaling coefficient of $\gamma_i = \alpha$ in \Cref{eq:scalednewton}. The scaled Newton polynomials defined by \Cref{eq:scalingcoeff}, on the contrary, take into account the relative position of each eigenvalue with the mean, leading to better scaling behaviors and, thus, larger maximum allowable step sizes.

In addition to the normalization effect seen in scalar scaling, column scaling further modifies the spectrum and results in a different distribution of eigenvalues. As a result, we see a significant change in block size for both Newton polynomials in the bottom two sub-figures in \Cref{fig:E20R5000_scaling}. In this case, the redistributed eigenvalues lead to a better quality of interpolation for Newton polynomials and give rise to increased step sizes. Monomials are only benefited by the normalization effect and therefore same level of improvement is observed for both matrix scaling methods. 

In summary, matrix equilibration generally has two combined effects on the eigenvalue spectrum. It not only attempts to normalize the spectrum to unit magnitude, but also modifies the relative eigenvalue distribution. We demonstrated the robustness of the scaled Newton polynomial, whose performance is insensitive to scalar scaling and depends mainly on the eigenvalue spectrum. This gives rise to a revised interpretation of the need for matrix equilibration that was deemed necessary in \cite{hoemmen2010communication,ballard2014communication,carson2015communication} for large step sizes. With the utilization of scaled Newton polynomials, attention can be directed towards the modification of eigenvalue distribution for large step sizes, which may be achieved through many alternative preconditioning methods.

\begin{figure}[tbhp]
\centering
\subfloat{
    \begin{tikzpicture}
    \begin{groupplot}[
        group style={
            group name=my fancy plots,
            group size=1 by 2,
            xticklabels at=edge bottom,
            vertical sep=0pt,
        },
        width=0.7\textwidth,
        xmin=0, xmax=600,
    ]
    \nextgroupplot[
        ymode=log,
        ymin=0.92,ymax=1,
        ytick={0.95,1},
        axis x line=top, 
        axis y discontinuity=parallel,
        height=0.25\textwidth,
        legend entries={GMRES,Monomial,Newton,Scaled Newton},
        legend columns=-1,
        legend style={at={(0.42,1.5)},anchor=north},
    ]
    \addplot [mark=o,mark repeat=40] table [x=i,y=GMRES] {Figures_data/E20R5000_Scalar/E20R5000Scalar_res.dat};
    \addplot [
    hcred,
    mark=diamond,
    mark repeat=40,
    mark phase=11,
    ] table [x=i,y=Monomial] {Figures_data/E20R5000_Scalar/E20R5000Scalar_res.dat};   
    \addplot [
    hcblue,
    mark=triangle,
    mark repeat=40,
    mark phase=21,
    ] table [x=i,y=Newton] {Figures_data/E20R5000_Scalar/E20R5000Scalar_res.dat};
    \addplot [
    hcgreen,
    mark=square,
    mark repeat=40,
    mark phase=31,
    ] table [x=i,y=SNewton] {Figures_data/E20R5000_Scalar/E20R5000Scalar_res.dat};
    
    \nextgroupplot[
        ymode=log,
        ymin=1e-16,ymax=1e-13,
        ytick={1e-16,1e-14},
        xtick={150,300,450,600},
        axis x line=bottom,
        height=0.25\textwidth,
        xlabel={iteration, $i$},
        ylabel={Relative residual / LOO},
        ylabel style={at={(ticklabel cs:1)}},
    ]
    \addplot [mark=o,mark repeat=40] table [x=i,y=GMRES] {Figures_data/E20R5000_Scalar/E20R5000Scalar_LOO.dat};
    \addplot [
    hcred, 
    mark=diamond,
    mark repeat=40,
    mark phase=11,
    ] table [x=i,y=Monomial] {Figures_data/E20R5000_Scalar/E20R5000Scalar_LOO.dat};     
    \addplot [
    hcblue, 
    mark=triangle,
    mark repeat=40,
    mark phase=21,
    ] table [x=i,y=Newton] {Figures_data/E20R5000_Scalar/E20R5000Scalar_LOO.dat};  
    \addplot [
    hcgreen, 
    mark=square,
    mark repeat=40,
    mark phase=31,
    ] table [x=i,y=SNewton] {Figures_data/E20R5000_Scalar/E20R5000Scalar_LOO.dat};  
    \end{groupplot}
    \end{tikzpicture}
}
\\\,\,
\subfloat{
    \begin{tikzpicture}
    \begin{axis}[
    width=0.7\textwidth,
    height=0.3\textwidth,
    xmin=1,
    xmax=600,
    xtick={150,300,450,600},
    ymin=0,
    ymax=120,
    ytick={0, 40, 80, 120},
    xlabel={iteration, $i$}, 
    ylabel={block size, $s$},
    ylabel shift=5,
    ]
    \addplot [
    hcred, 
    mark=diamond,
    mark repeat=30,
    ] table [x=i,y=Monomial] {Figures_data/E20R5000_Scalar/E20R5000Scalar_s.dat};
    \addplot [
    hcblue, 
    mark=triangle,
    mark repeat=30,
    ] table [x=i,y=Newton] {Figures_data/E20R5000_Scalar/E20R5000Scalar_s.dat};
    \addplot [
    hcgreen, 
    mark=square,
    mark repeat=30,
    ] table [x=i,y=SNewton] {Figures_data/E20R5000_Scalar/E20R5000Scalar_s.dat};
    \end{axis}
    \end{tikzpicture}
}
\\
\subfloat{
    \begin{tikzpicture}
    \begin{groupplot}[
        group style={
            group name=my fancy plots,
            group size=1 by 2,
            xticklabels at=edge bottom,
            vertical sep=0pt,
        },
        width=0.7\textwidth,
        xmin=0, xmax=600,
    ]
    \nextgroupplot[
        ymode=log,
        ymin=0.92,ymax=1,
        ytick={0.95,1},
        axis x line=top, 
        axis y discontinuity=parallel,
        height=0.25\textwidth,
    ]
    \addplot [mark=o, mark repeat=40] table [x=i,y=GMRES] {Figures_data/E20R5000_Column/E20R5000Column_res.dat};
    \addplot [
    hcred,
    mark=diamond,
    mark repeat=40,
    mark phase=11,
    ] table [x=i,y=Monomial] {Figures_data/E20R5000_Column/E20R5000Column_res.dat};   
    \addplot [
    hcblue,
    mark=triangle,
    mark repeat=40,
    mark phase=21,
    ] table [x=i,y=Newton] {Figures_data/E20R5000_Column/E20R5000Column_res.dat};
    \addplot [
    hcgreen,
    mark=square,
    mark repeat=40,
    mark phase=31,
    ] table [x=i,y=SNewton] {Figures_data/E20R5000_Column/E20R5000Column_res.dat};
    
    \nextgroupplot[
        ymode=log,
        ymin=1e-16,ymax=1e-13,
        ytick={1e-16,1e-14},
        xtick={150,300,450,600},
        axis x line=bottom,
        height=0.25\textwidth,
        xlabel={iteration, $i$},
        ylabel={Relative residual / LOO},
        ylabel style={at={(ticklabel cs:1)}},
    ]
    \addplot [mark=o,mark repeat=40] table [x=i,y=GMRES] {Figures_data/E20R5000_Column/E20R5000Column_LOO.dat};
    \addplot [
    hcred, 
    mark=diamond,
    mark repeat=40,
    mark phase=11,
    ] table [x=i,y=Monomial] {Figures_data/E20R5000_Column/E20R5000Column_LOO.dat};     
    \addplot [
    hcblue, 
    mark=triangle,
    mark repeat=40,
    mark phase=21,
    ] table [x=i,y=Newton] {Figures_data/E20R5000_Column/E20R5000Column_LOO.dat};  
    \addplot [
    hcgreen, 
    mark=square,
    mark repeat=40,
    mark phase=31,
    ] table [x=i,y=SNewton] {Figures_data/E20R5000_Column/E20R5000Column_LOO.dat};  
    \end{groupplot}
    \end{tikzpicture}
}
\\\,\,
\subfloat{
    \begin{tikzpicture}
    \begin{axis}[
    width=0.7\textwidth,
    height=0.3\textwidth,
    xmin=1,
    xmax=600,
    xtick={150,300,450,600},
    ymin=0,
    ymax=150,
    ytick={0, 50, 100, 150},
    xlabel={iteration, $i$}, 
    ylabel={block size, $s$},
    ylabel shift=5,
    ]
    \addplot [
    hcred, 
    mark=diamond,
    mark repeat=30,
    ] table [x=i,y=Monomial] {Figures_data/E20R5000_Column/E20R5000Column_s.dat};
    \addplot [
    hcblue,
    mark=triangle,
    mark repeat=30,
    ] table [x=i,y=Newton] {Figures_data/E20R5000_Column/E20R5000Column_s.dat};
    \addplot [
    hcgreen,
    mark=square,
    mark repeat=30,
    ] table [x=i,y=SNewton] {Figures_data/E20R5000_Column/E20R5000Column_s.dat};
    \end{axis}
    \end{tikzpicture}
}
\caption{E20R5000 matrix with scalar scaling (top 2 sub-figures) and column scaling (bottom 2 sub-figures). The upper sub-figure in each case shows the relative residual and LOO errors of GMRES and adaptive $s$-step GMRES using different polynomial bases in the MPK. The lower sub-figure shows the adapted step size $s$.}
\label{fig:E20R5000_scaling}
\end{figure}
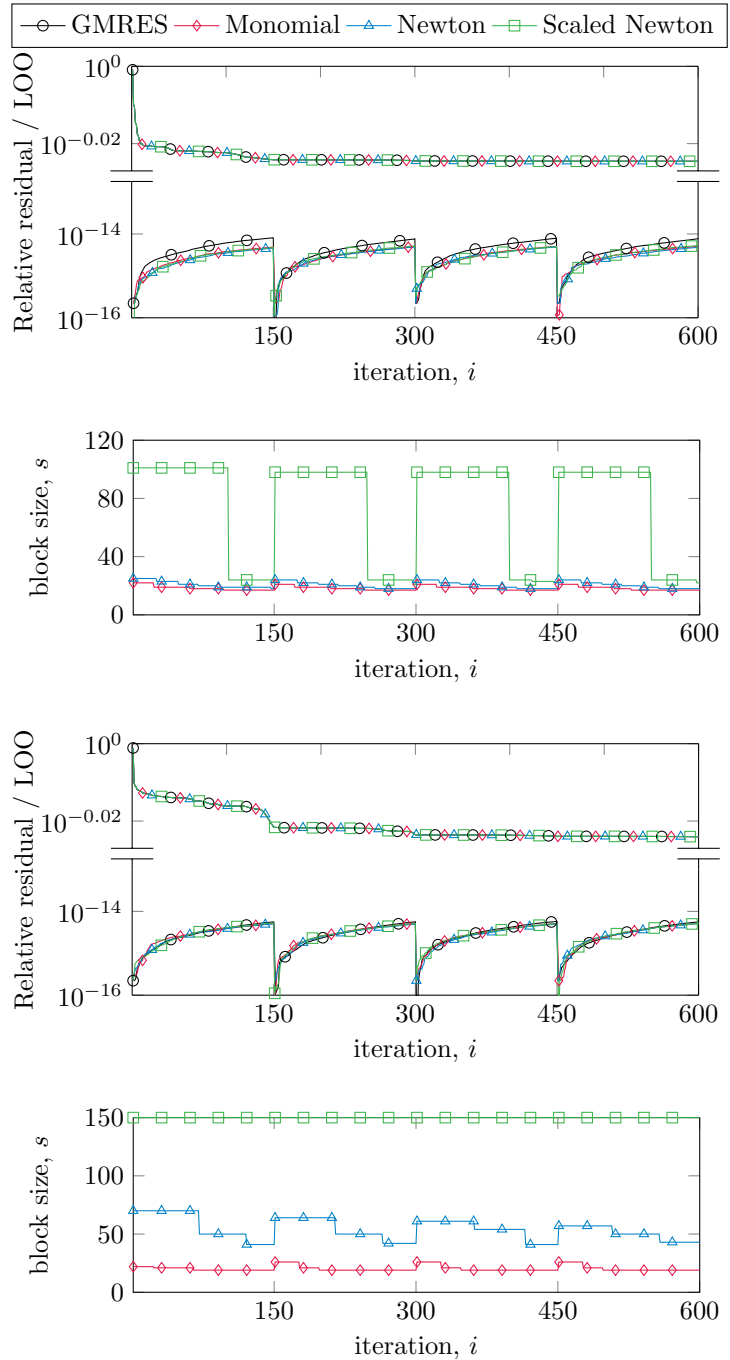 %% lengthy figures are isolated into individual files
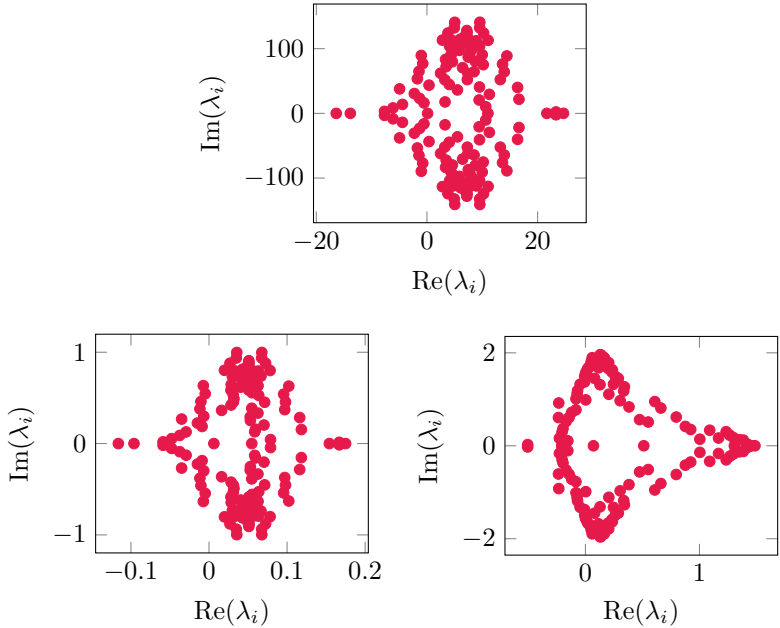
\begin{figure}[tbhp]
\centering
\subfloat{
    \begin{tikzpicture}
    \begin{axis}[
    width=0.4\textwidth,
    xlabel={Re$(\lambda_i)$}, 
    ylabel={Im$(\lambda_i)$},
    ]
    \addplot [
    only marks,
    hcred,
    mark=*,
    % mark repeat=5,
    ] table [x=eig_real,y=eig_imag] {Figures_data/E20R5000/E20R5000_spectrum.dat};
    \end{axis}
    \end{tikzpicture}
}
\\
\subfloat{
    \begin{tikzpicture}
    \begin{axis}[
    width=0.4\textwidth,
    xlabel={Re$(\lambda_i)$}, 
    ylabel={Im$(\lambda_i)$},
    ]
    \addplot [
    only marks,
    hcred,
    mark=*,
    % mark repeat=5,
    ] table [x=eig_real,y=eig_imag] {Figures_data/E20R5000_Scalar/E20R5000Scalar_spectrum.dat};
    \end{axis}
    \end{tikzpicture}
}
\subfloat{
    \begin{tikzpicture}
    \begin{axis}[
    width=0.4\textwidth,
    xlabel={Re$(\lambda_i)$}, 
    ylabel={Im$(\lambda_i)$},
    ]
    \addplot [
    only marks,
    hcred,
    mark=*,
    % mark repeat=5,
    ] table [x=eig_real,y=eig_imag] {Figures_data/E20R5000_Column/E20R5000Column_spectrum.dat};
    \end{axis}
    \end{tikzpicture}
}
\caption{Eigenvalue spectrum of matrix E20R5000. The three sub-figures show the eigenvalue spectrum of the same matrix but under different preconditioning techniques. Top: no preconditioning. Lower left: scalar scaling. Lower right: column scaling.}
\label{fig:E20R5000_spectrum}
\end{figure} %% lengthy figures are isolated into individual files

\subsubsection{2D Laplace matrix with ILU preconditioner}
As most test cases presented so far were not preconditioned with an effective preconditioner, they have not been converged down to sufficiently small residuals (except for tests with the diagonal matrix). For completeness, we used subsequent tests to demonstrate the stability of \Cref{alg:adaptsgmres} up to a satisfactory level of convergence.

\begin{figure}[tbhp]
\centering
\subfloat{
    \begin{tikzpicture}
    \begin{axis}[
    width=0.7\textwidth,
    ymode=log,
    xmin=1,
    xmax=400,
    ymin=1e-16,
    ymax=1,
    xlabel={iteration, $i$}, 
    ylabel={Relative residual / LOO},
    legend entries={GMRES, Monomial, Newton, Scaled Newton},
    legend columns=-1,
    legend style={at={(0.5,1.15)},anchor=north},
    ]
    \addplot [mark=o,mark repeat=24] table [x=i,y=GMRES] {Figures_data/2DLaplace_ILU/2DLaplaceILU_res.dat};
    \addplot [
    hcred,
    mark=diamond,
    mark repeat=24,
    mark phase=7,
    ] table [x=i,y=Monomial] {Figures_data/2DLaplace_ILU/2DLaplaceILU_res.dat};
    \addplot [
    hcblue,
    mark=triangle,
    mark repeat=24,
    mark phase=13,
    ] table [x=i,y=Newton] {Figures_data/2DLaplace_ILU/2DLaplaceILU_res.dat};
    \addplot [
    hcgreen, 
    mark=square,
    mark repeat=24,
    mark phase=19,
    ] table [x=i,y=SNewton] {Figures_data/2DLaplace_ILU/2DLaplaceILU_res.dat};
    % \addplot [
    % hcgrey, dashed,
    % mark=*,
    % mark repeat=24,
    % ] table [x=i,y=Chebyshev] {Figures_data/2DLaplace_ILU/2DLaplaceILU_res.dat};
    \addplot [mark=o,mark repeat=24] table [x=i,y=GMRES] {Figures_data/2DLaplace_ILU/2DLaplaceILU_LOO.dat};
    \addplot [
    hcred, 
    mark=diamond,
    mark repeat=24,
    mark phase=7,
    ] table [x=i,y=Monomial] {Figures_data/2DLaplace_ILU/2DLaplaceILU_LOO.dat};
    \addplot [
    hcblue, 
    mark=triangle,
    mark repeat=24,
    mark phase=13,
    ] table [x=i,y=Newton] {Figures_data/2DLaplace_ILU/2DLaplaceILU_LOO.dat};
    \addplot [
    hcgreen,
    mark=square,
    mark repeat=24,
    mark phase=19,
    ] table [x=i,y=SNewton] {Figures_data/2DLaplace_ILU/2DLaplaceILU_LOO.dat};
    % \addplot [
    % hcgrey, dashed,
    % mark=*,
    % mark repeat=24,
    % ] table [x=i,y=Chebyshev] {Figures_data/2DLaplace_ILU/2DLaplaceILU_LOO.dat};
    \end{axis}
    \end{tikzpicture}
}
\\
\subfloat{
    \begin{tikzpicture}
    \begin{groupplot}[
        group style={
            group name=my fancy plots,
            group size=1 by 2,
            xticklabels at=edge bottom,
            vertical sep=0pt,
        },
        width=0.7\textwidth,
        xmin=0, xmax=400,
    ]
    \nextgroupplot[
        ymin=375,ymax=400,
        axis x line=top, 
        axis y discontinuity=parallel,
        height=0.2\textwidth,
        ytick={390,400},
    ]
    \addplot [
    hcgreen,
    mark=square,
    mark repeat=20,
    ] table [x=i,y=SNewton] {Figures_data/2DLaplace_ILU/2DLaplaceILU_s.dat};     
    
    \nextgroupplot[
        ymin=0,ymax=25,
        axis x line=bottom,
        x axis line style={},
        height=0.2\textwidth,
        xlabel={iteration, $i$},
        ylabel={block size, $s$},
        ylabel style={at={(ticklabel cs:1)}},
        ytick={0,10,20},
        ylabel shift=6,
    ]
    \addplot [
    hcred,
    mark=diamond,
    mark repeat=20,
    ] table [x=i,y=Monomial] {Figures_data/2DLaplace_ILU/2DLaplaceILU_s.dat};   
    \addplot [
    hcblue,
    mark=triangle,
    mark repeat=20,
    ] table [x=i,y=Newton] {Figures_data/2DLaplace_ILU/2DLaplaceILU_s.dat};
    % \addplot [
    % hcgrey, dashed,
    % mark=*,
    % mark repeat=20,
    % ] table [x=i,y=Chebyshev] {Figures_data/2DLaplace_ILU/2DLaplaceILU_s.dat};
    \end{groupplot}
    \end{tikzpicture}
}
\caption{2D Laplace matrix with ILU(0) preconditioner. Top figure: relative residual and LOO errors of GMRES and adaptive $s$-step GMRES with different bases. Convergence is accelerated by preconditioning and LOO error of the adaptive s-step GMRES is still kept near machine epsilon. Lower figure: adapted step size, $s$. It is worth mentioning that scaled Newton polynomials give a block size of $s=400$, which is equivalent to one block iteration for convergence.}
\label{fig:test_laplace_precond}
\end{figure}
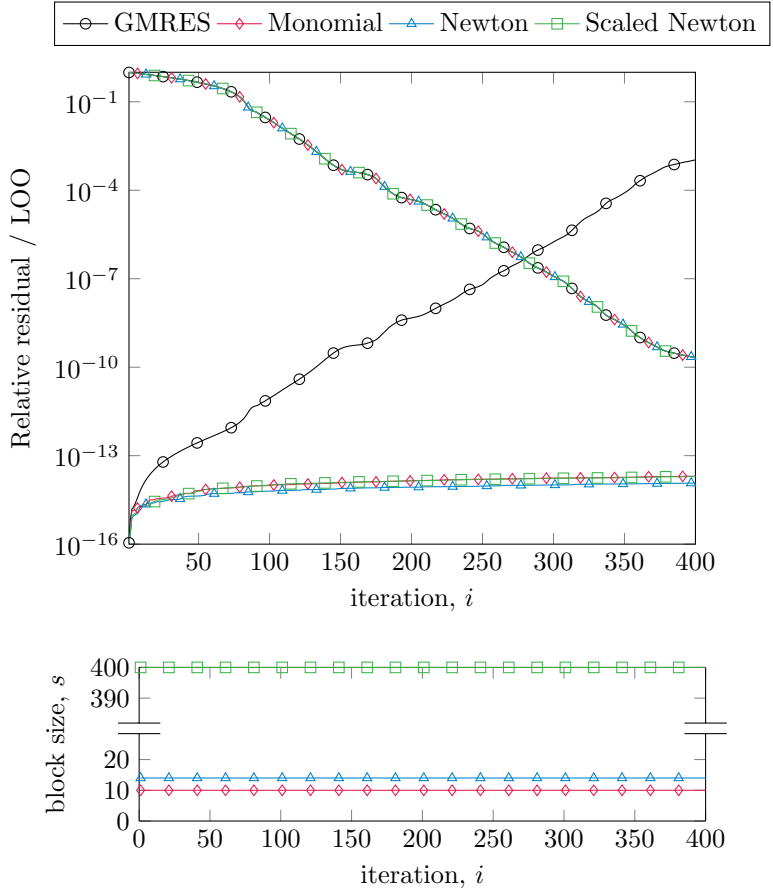

The same matrix problem as \Cref{sec:test_2Dlaplace} is used with ILU(0) preconditioner. The initial step size is $s_0 = 400$. As shown in \Cref{fig:test_laplace_precond}, convergence is accelerated by the presence of a preconditioner. However, all implementations with various polynomial bases agree with the GMRES baseline. A minimal LOO error is observed for all cases. It is worth mentioning that the adapted step size for scaled Newton polynomials is 400, indicating that one single block iteration is sufficient for the linear system to converge with guaranteed stability.

 %% lengthy figures are isolated into individual files

\subsubsection{E20R5000 matrix with ILU preconditioner}\label{sec:test_E20R5000_ILU}
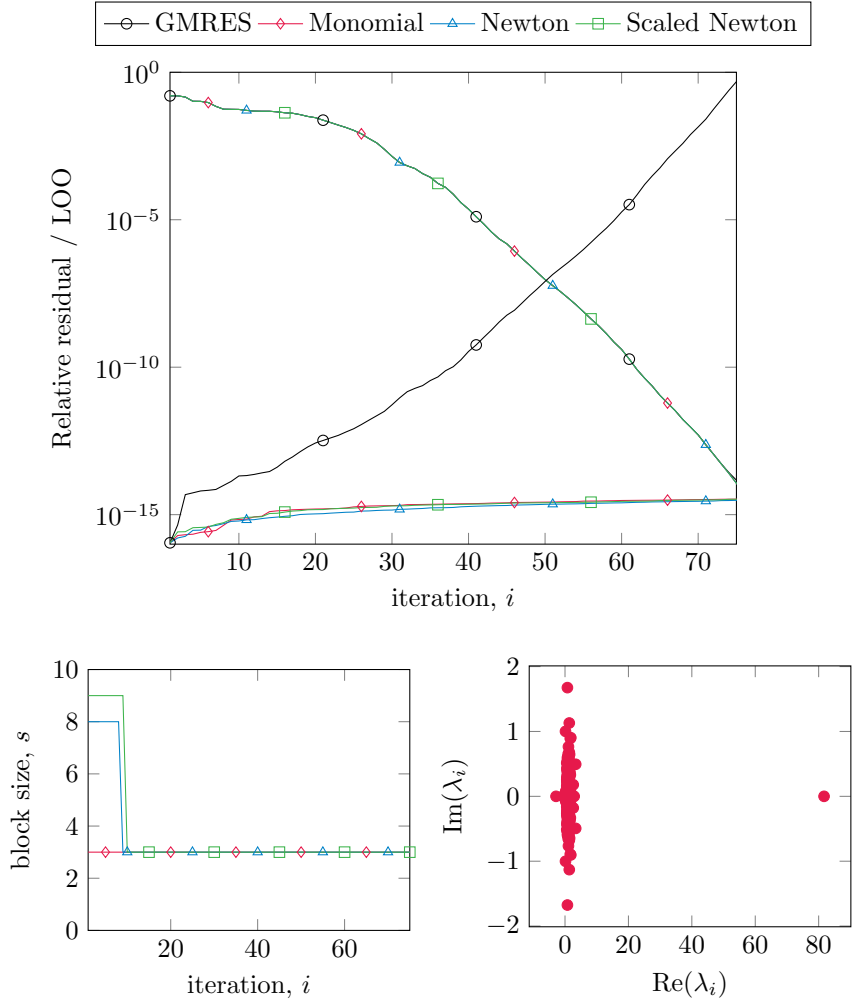
\begin{figure}[tbhp]
\centering
\subfloat{
    \begin{tikzpicture}
    \begin{axis}[
    width=0.7\textwidth,
    ymode=log,
    xmin=1,
    xmax=75,
    ytick={1e-15, 1e-10, 1e-5, 1},
    ymax=1,
    ymin=1e-16,
    xlabel={iteration, $i$}, 
    ylabel={Relative residual / LOO},
    legend entries={GMRES,Monomial,Newton,Scaled Newton},
    legend columns=-1,
    legend style={at={(0.5,1.15)},anchor=north},
    % legend pos=north east
    ]
    \addplot [mark=o,mark repeat=20] table [x=i,y=GMRES] {Figures_data/E20R5000_ILU/E20R5000ILU_res.dat};
    \addplot [
    hcred, 
    mark=diamond,
    mark repeat=20,
    mark phase=6,
    ] table [x=i,y=Monomial] {Figures_data/E20R5000_ILU/E20R5000ILU_res.dat};
    \addplot [
    hcblue, 
    mark=triangle,
    mark repeat=20,
    mark phase=11,
    ] table [x=i,y=Newton] {Figures_data/E20R5000_ILU/E20R5000ILU_res.dat};
    \addplot [
    hcgreen, 
    mark=square,
    mark repeat=20,
    mark phase=16,
    ] table [x=i,y=SNewton] {Figures_data/E20R5000_ILU/E20R5000ILU_res.dat};
    \addplot [mark=o,mark repeat=20] table [x=i,y=GMRES] {Figures_data/E20R5000_ILU/E20R5000ILU_LOO.dat};
    \addplot [
    hcred, 
    mark=diamond,
    mark repeat=20,
    mark phase=6,
    ] table [x=i,y=Monomial] {Figures_data/E20R5000_ILU/E20R5000ILU_LOO.dat};
    \addplot [
    hcblue, 
    mark=triangle,
    mark repeat=20,
    mark phase=11,
    ] table [x=i,y=Newton] {Figures_data/E20R5000_ILU/E20R5000ILU_LOO.dat};
    \addplot [
    hcgreen, 
    mark=square,
    mark repeat=20,
    mark phase=16,
    ] table [x=i,y=SNewton] {Figures_data/E20R5000_ILU/E20R5000ILU_LOO.dat};
    \end{axis}
    \end{tikzpicture}
}\\
\subfloat{
    \begin{tikzpicture}
    \begin{axis}[
    width=0.45\textwidth,
    % height=0.3\textwidth,
    xmin=1,
    xmax=75,
    ymin=0,
    ymax=10,
    xlabel={iteration, $i$}, 
    ylabel={block size, $s$},
    % ylabel shift=5,
    ]
    \addplot [
    hcred, 
    mark=diamond,
    mark repeat=15,
    mark phase=5,
    ] table [x=i,y=Monomial] {Figures_data/E20R5000_ILU/E20R5000ILU_s.dat};
    \addplot [
    hcblue, 
    mark=triangle,
    mark repeat=15,
    mark phase=10,
    ] table [x=i,y=Newton] {Figures_data/E20R5000_ILU/E20R5000ILU_s.dat};
    \addplot [
    hcgreen, 
    mark=square,
    mark repeat=15,
    mark phase=15,
    ] table [x=i,y=SNewton] {Figures_data/E20R5000_ILU/E20R5000ILU_s.dat};
    % \addplot [
    % hcgrey, dashed,
    % mark=*,
    % mark repeat=7,
    % ] table [x=i,y=Chebyshev] {Figures_data/E20R5000_ILU/E20R5000ILU_s.dat};
    \end{axis}
    \end{tikzpicture}
}
\subfloat{
    \begin{tikzpicture}
    \begin{axis}[
    width=0.45\textwidth,
    xlabel={Re$(\lambda_i)$}, 
    ylabel={Im$(\lambda_i)$},
    ]
    \addplot [
    only marks,
    hcred,
    mark=*,
    % mark repeat=5,
    ] table [x=eig_real,y=eig_imag] {Figures_data/E20R5000_ILU/E20R5000ILU_spectrum.dat};
    \end{axis}
    \end{tikzpicture}
}
\caption{E20R5000 matrix with the ILUTP preconditioner. Top figure: relative residual and LOO errors. A rapid acceleration of the convergence is observed due to the effective preconditioning, while the stability of the orthogonalization is still guaranteed. Lower left figure: adapted step size, $s$. Adapted block sizes are generally very small. Lower right figure: eigenvalue spectrum of the preconditioned matrix. One outlying eigenvalue is well-separated from the rest of the eigenvalues.}
\label{fig:E20R5000_ILU}
\end{figure}

A linear system involving the E20R5000 matrix with the ILU(0) preconditioner does not converge well. As an alternative, we present the results in which the matrix is preconditioned using ILUTP (ILU with threshold and pivoting) with a drop tolerance of $10^{-4}$. Using $s_0 = 75$, an accelerated convergence is observed in \Cref{fig:E20R5000_ILU} with stable orthogonalization. Even though the matrix is better conditioned for convergence, the modified eigenvalue spectrum in \Cref{fig:E20R5000_ILU} shows an outlying eigenvalue around $\lambda_\text{max} = 80$ that is separated from the rest of the spectrum. Such outlying eigenvalues degrade the quality of eigenvalues as interpolation nodes, and the incremental condition number grows rapidly. Hence, a much smaller step size is observed, even for scaled Newton polynomials. This shows that preconditioning does not necessarily lead to larger step sizes for polynomial bases in the MPK. 

 %% lengthy figures are isolated into individual files

\subsection{Numerical experiments on initial step size estimator}\label{subsec:tests_estimator}

Here, we incorporate the scaled Newton polynomial in \Cref{subsec:scalednewton} with the initial step size estimator in \Cref{subsec:errorestimator} for the adaptive $s$-step GMRES algorithm in \Cref{alg:adaptsgmres}. We demonstrate the initial step size estimator using four detailed examples, and the results are presented in \Cref{fig:errorestimator}. In all four examples, we show the eigenvalue spectrum on the left. On the right, we plot the incremental condition number of the Krylov basis matrix $V$ generated by scaled Newton polynomials (i.e., $\kappa(V_{1:j})$), the column-wise norm of $V$ (i.e., $\|V_j\|$), and the corresponding column-wise norm prediction $\|E_j\|$ from the initial step size estimator.

\subsubsection{Two diagonal matrices with known eigenvalues}
We first use two examples of diagonal matrices with known eigenvalues. Both diagonal matrices are of size $N=200$, and we use $s_0 = 200$ in both tests. The first matrix has uniform eigenvalues $\lambda_i \in (1, 200)$. The second matrix has the same eigenvalue distribution except that the largest eigenvalue is increased to $\lambda_\text{max} = 2000$. In both test cases, the estimator prediction $\|E_j\|$ shows the same trend as the actual column-wise norm $\|V_j\|$. As an increasing column-wise norm inevitably increases the incremental condition number, $\kappa(V_{1:j})$ grows at a faster rate than the predicted and actual column-wise norms. The incremental condition number of the first test matrix increases beyond $O(10^8)$ at $j = 62$, indicating an initial step size of $s_0 = 200$ would have to be truncated down to 62 Krylov basis vectors for stable orthogonalization. The initial step size estimator based on column-wise norms predicts $s_0^* = 134$, significantly reducing the number of wasted Krylov vectors. Similarly, the initial step size estimator for the second matrix predicts $s_0^* = 17$ while $\kappa(V_{1:j})$ reaches $O(10^8)$ at $j=15$.

To illustrate the necessity of considering numerical errors in the estimator, we also include $\|E^L_j\|$ in the two test cases in \Cref{fig:errorestimator}, where $E^L$ is the strictly lower triangular part of the auxiliary matrix $E$. In other words, $E^L$ denotes the evolution of the product term when Ritz values are exactly equal to the true eigenvalues without roundoff errors. In both tests, $\|E^L_j\|$, unlike $\|E_j\|$, stays at a relatively constant order of magnitude without growing exponentially. This indicates that the exponential growth of $\|V_j\|$ is largely due to roundoff errors in Ritz value approximations. The onset of exponential growth is only captured if we consider the evolution of roundoff errors in the upper triangular part of $E$, justifying the definition of the auxiliary matrix $E$ in~\eqref{eq:auxiliarymatrix}. The exact column at which exponential growth starts depends on the eigenvalue spectrum. In the second test case with $\lambda_\text{max} = 2000$, the roundoff error at $i = 1$ is constantly increased by a factor of $O(10)$ with every $j$, leading to the rapid growth of roundoff error and therefore small initial step size.

\subsubsection{Unpreconditioned and preconditioned E20R5000 matrices}
\begin{figure}[!tbhp]
\centering
\begin{tikzpicture}
\begin{groupplot}[
    group style={
        group size=2 by 4,
        horizontal sep=1.5cm,
        vertical sep=1.5cm,
        },
]
\nextgroupplot[
    width=0.4\textwidth,
    % height=0.3\textwidth,
    ymin=-1,
    ymax=1,
    xlabel={Re$(\lambda_i)$}, 
    ylabel={Im$(\lambda_i)$},
    ]
    \addplot [
    hcred, only marks,
    mark=*,
    mark repeat=10,
    ] table [x=eig_real,y=eig_imag] {Figures_data/ErrorEstimator/Diagonal_200_error_updated.dat};
    
\nextgroupplot[
    ymode=log,
    width=0.4\textwidth,
    % height=0.3\textwidth,
    xmin=1,
    xmax=150,
    ymin=1e-2,
    ymax=1e16,
    xlabel={column, $j$}, 
    ytick={1, 1e8, 1e16},
    ]
    \addplot [
    hcgreen, 
    mark=*,
    mark repeat=12,
    ] table [x=i,y=cond] {Figures_data/ErrorEstimator/Diagonal_200_error_updated.dat};
    \addplot [
    hcblue,
    mark=*,
    mark repeat=12,
    ] table [x=i,y=vecnorm] {Figures_data/ErrorEstimator/Diagonal_200_error_updated.dat};
    \addplot [
    black, 
    mark=triangle*,
    mark repeat=12,
    ] table [x=i,y=E] {Figures_data/ErrorEstimator/Diagonal_200_error_updated.dat};
    \addplot [
    hcorange,
    mark=triangle*,
    mark repeat=12,
    ] table [x=i,y=E_true] {Figures_data/ErrorEstimator/Diagonal_200_error_updated.dat};
    
\nextgroupplot[
    width=0.4\textwidth,
    % height=0.3\textwidth,
    ymin=-1,
    ymax=1,
    xlabel={Re$(\lambda_i)$}, 
    ylabel={Im$(\lambda_i)$},
    ]
    \addplot [
    hcred, only marks,
    mark=*,
    % mark repeat=10,
    ] table [x=eig_real,y=eig_imag] {Figures_data/ErrorEstimator/Diagonal_2000_error_updated.dat};
    
\nextgroupplot[
    ymode=log,
    width=0.4\textwidth,
    % height=0.3\textwidth,
    xmin=1,
    xmax=30,
    ymin=1e-2,
    ymax=1e16,
    xlabel={column, $j$}, 
    ytick={1, 1e8, 1e16},
    legend entries={$\kappa(V_{:,1:j})$,$\|V_{:,j}\|$,$\|E_j\|$,$\|E^L_j\|$},
    legend style={at={(1.7,-0.2)},anchor=east},
    ]
    \addplot [
    hcgreen, 
    mark=*,
    mark repeat=2,
    ] table [x=i,y=cond] {Figures_data/ErrorEstimator/Diagonal_2000_error_updated.dat};
    \addplot [
    hcblue, 
    mark=*,
    mark repeat=2,
    ] table [x=i,y=vecnorm] {Figures_data/ErrorEstimator/Diagonal_2000_error_updated.dat};
    \addplot [
    black, 
    mark=triangle*,
    mark repeat=2,
    ] table [x=i,y=E] {Figures_data/ErrorEstimator/Diagonal_2000_error_updated.dat};
    \addplot [
    hcorange, 
    mark=triangle*,
    mark repeat=2,
    ] table [x=i,y=E_true] {Figures_data/ErrorEstimator/Diagonal_2000_error_updated.dat};
    
\nextgroupplot[
    width=0.4\textwidth,
    % height=0.3\textwidth,
    xlabel={Re$(\lambda_i)$}, 
    ylabel={Im$(\lambda_i)$},
    ]
    \addplot [
    hcred, only marks,
    mark=*,
    % mark repeat=10,
    ] table [x=eig_real,y=eig_imag] {Figures_data/E20R5000/E20R5000_spectrum.dat};
    
\nextgroupplot[ymode=log,
    width=0.4\textwidth,
    % height=0.3\textwidth,
    xmin=1,
    xmax=150,
    ymin=1e-2,
    ymax=1e16,
    xlabel={column, $j$}, 
    ytick={1, 1e8, 1e16},
    ]
    \addplot [
    hcgreen, 
    mark=*,
    mark repeat=12,
    ] table [x=i,y=cond] {Figures_data/ErrorEstimator/E20R5000_error_updated.dat};
        \addplot [
    hcblue, 
    mark=*,
    mark repeat=12,
    ] table [x=i,y=vecnorm] {Figures_data/ErrorEstimator/E20R5000_error_updated.dat};
    \addplot [
    black, 
    mark=triangle*,
    mark repeat=12,
    ] table [x=i,y=E] {Figures_data/ErrorEstimator/E20R5000_error_updated.dat};

\nextgroupplot[
    width=0.4\textwidth,
    % height=0.3\textwidth,
    xlabel={Re$(\lambda_i)$}, 
    ylabel={Im$(\lambda_i)$},
    ]
    \addplot [
    hcred, only marks,
    mark=*,
    % mark repeat=10,
    ] table [x=eig_real,y=eig_imag] {Figures_data/E20R5000_ILU/E20R5000ILU_spectrum.dat};

\nextgroupplot[
ymode=log,
    width=0.4\textwidth,
    % height=0.3\textwidth,
    xmin=1,
    xmax=20,
    ymin=1e-2,
    ymax=1e16,
    xlabel={column, $j$}, 
    xtick={2,4,6,8,10,12,14,16,18},
    ytick={1, 1e8, 1e16},
    ]
        \addplot [
    hcgreen, 
    mark=*,
    % mark repeat=6,
    ] table [x=i,y=cond] {Figures_data/ErrorEstimator/E20R5000ILU_error_updated.dat};
    \addplot [
    hcblue, 
    mark=*,
    % mark repeat=6,
    ] table [x=i,y=vecnorm] {Figures_data/ErrorEstimator/E20R5000ILU_error_updated.dat};
    \addplot [
    black, 
    mark=triangle*,
    % mark repeat=6,
    ] table [x=i,y=E] {Figures_data/ErrorEstimator/E20R5000ILU_error_updated.dat};

\end{groupplot}
\end{tikzpicture} 
\caption{Eigenvalue spectra (left) and performance of initial step size estimators (right) for four matrix problems (from top to bottom): (1) diagonal matrix with eigenvalues uniformly in (1, 200); (2) same diagonal matrix as the previous one but the largest eigenvalue is increased to 2,000; (3) E20R5000; (4) E20R5000 with ILUTP preconditioner.}
\label{fig:errorestimator}
\end{figure}
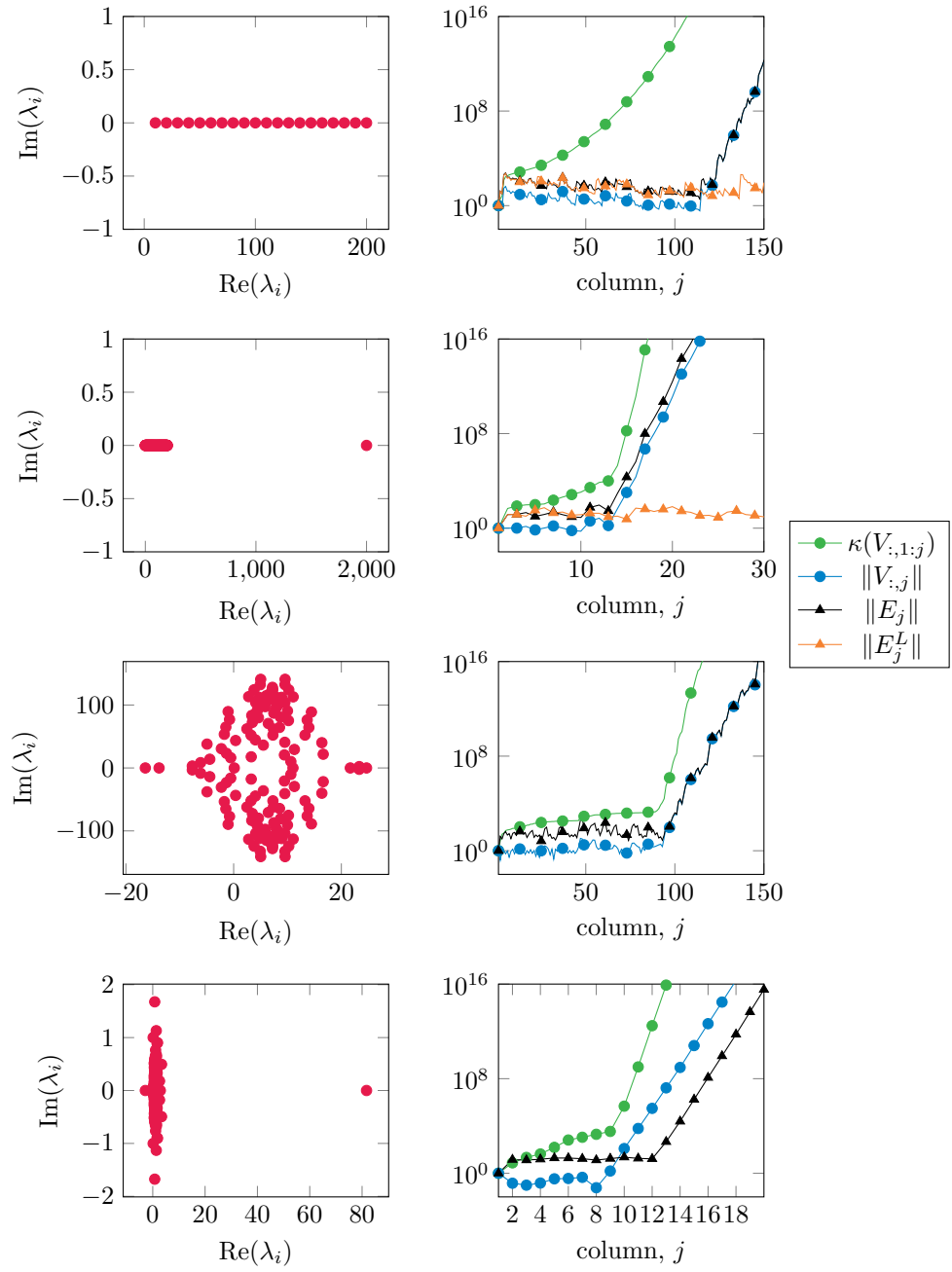

The next two examples are the unpreconditioned E20R5000 matrix from \Cref{sec:test_E20R5000} and ILUTP-preconditioned E20R5000 matrix from \Cref{sec:test_E20R5000_ILU}. The eigenvalue spectra are reproduced in \Cref{fig:errorestimator} for ease of comparison. The predictions of the numerical error given by the initial step size estimator in both cases show excellent agreement with the growth rate of the column-wise vector norm of the Krylov basis vectors. A threshold of $\Omega_\text{est}=10^7$ restricts the initial step size to $s_0^*=113$ and $s_0^*=13$ for the two cases, respectively. Both estimates are very close to the adapted step size shown in \Cref{fig:E20R5000,fig:E20R5000_ILU}.

In all examples, we show that the initial step size estimator correctly captures the growth rate of the column-wise norm of Krylov basis matrices. Since the increased vector norm unavoidably increases the incremental condition number, one could use the estimator to upper-bound the initial step size such that the vector norm does not increase beyond a certain threshold. This, in turn, limits the growth of the incremental condition number and leads to a well-informed step size, ultimately reducing the number of unnecessarily generated Krylov basis vectors.

\subsubsection{Matrices from the SuiteSparse Matrix Collection}\label{sec:tests_final}

\begin{figure}
    \centering
    \begin{tikzpicture}
    \begin{axis}[
        width=.9\textwidth,
        height=.5\textwidth,
        bar width=3pt,
        ybar,
        enlarge x limits=0.1,
        % xmin=cage14,
        % xmax=NACA0015,
        xmin=atmosmodd,
        xmax=thermal2,
        ymin=0,
        ymax=550,
        ytick={0,100,200,300,400,500},
        ylabel={block size, $s$},
        % symbolic x coords={cage14, Hamrles3, G3\_circuit, thermal2, ecology1, atmosmodd, atmosmodj, atmosmodl, atmosmodm, roadNet-CA, roadNet-PA, roadNet-TX, delaunay\_n20, belgium\_osm, NACA0015},
        symbolic x coords={atmosmodd, atmosmodj, atmosmodl, cage14, webbase-1M, G3\_circuit, Hamrle3, ecology1, af\_shell10, Transport, CurlCurl\_3, CurlCurl\_4, nlpkkt80, nlpkkt240, thermal2},
        xtick=data,
        xticklabel style={rotate=90},
        % nodes near coords,
        % nodes near coords align={vertical},
        ymajorgrids=true,
    ]
    \addplot [
    fill = black,
    % ] coordinates {(cage14,17) (Hamrles3,500) (G3\_circuit,14) (thermal2, 169) (ecology1, 246) (atmosmodd, 75) (atmosmodj, 76) (atmosmodl, 68) (atmosmodm, 58) (roadNet-CA, 122) (roadNet-PA, 139)(roadNet-TX, 106)  (delaunay\_n20, 116) (belgium\_osm, 128) (NACA0015, 313)}; % actual
    ] coordinates { (atmosmodd, 75) (atmosmodj, 76) (atmosmodl, 68) (cage14, 17) (webbase-1M, 66) (G3\_circuit, 14) (Hamrle3, 500) (ecology1, 246) (af\_shell10, 71) (Transport, 147) (CurlCurl\_3, 29) (CurlCurl\_4, 29) (nlpkkt80, 146) (nlpkkt240, 146) (thermal2, 172) }; %actual
    \addplot [
    pattern = dots,
    ] coordinates { (atmosmodd, 80) (atmosmodj, 84) (atmosmodl, 70) (cage14, 48) (webbase-1M, 148) (G3\_circuit, 35) (Hamrle3, 500) (ecology1, 257) (af\_shell10, 58) (Transport, 197) (CurlCurl\_3, 41) (CurlCurl\_4, 38) (nlpkkt80, 158) (nlpkkt240, 157) (thermal2, 190) }; %estimator
    % ] coordinates {(cage14,48) (Hamrles3,500) (G3\_circuit,35) (thermal2, 190) (ecology1, 257) (atmosmodd, 80) (atmosmodj, 85) (atmosmodl, 84) (atmosmodm, 70) (roadNet-CA, 138) (roadNet-PA, 154) (roadNet-TX, 120) (delaunay\_n20, 134) (belgium\_osm, 148) (NACA0015, 360) }; % estimator
    \legend{Adapted step size, Estimator prediction ($s_0^*$)}
    \end{axis}
    \end{tikzpicture}
    \caption{Matrices from SuiteSparse Matrix Collection. The adaptive $s$-step GMRES algorithm with scaled Newton polynomials is used in all cases. With an initial step size of $s_0=500$, the adapted step size as well as the predicted step size given by the estimator $s_0^*$ are shown for each matrix. The larger the block size, the more communication savings can be achieved. The smaller the discrepancy between two step sizes, the less computation and/or communication is wasted.}
    \label{fig:suitesparse}
\end{figure}
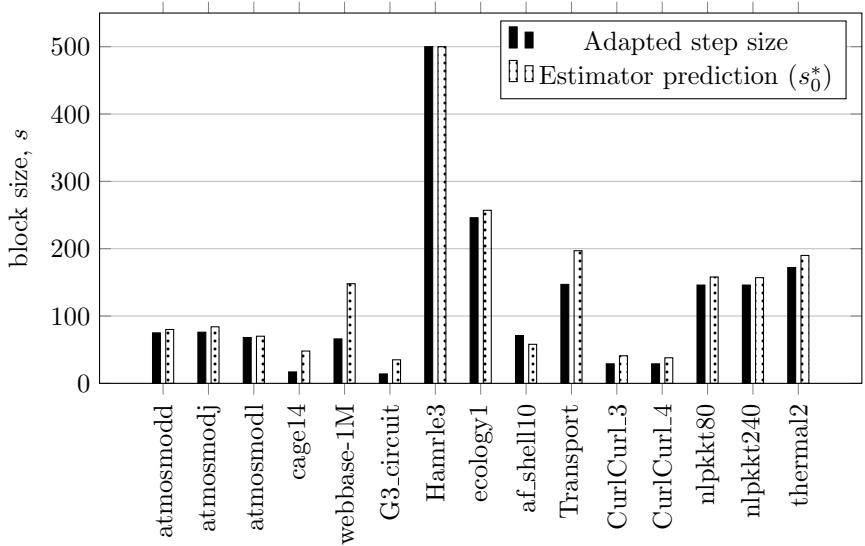

In this series of tests, we show the robustness of the overall algorithm by testing the framework on a wide range of matrices from the SuiteSparse Matrix Collection~\cite{davis2011university}. All matrices are of size $N>10^6$. Properties of the test matrices are shown in~\Cref{tab:matrix_properties}. Similar to the tests presented above, the relative residual agrees well with the GMRES algorithm with minimal LOO for all tests. With an initial step size of $s_0=500$, we plot the maximum adapted step size for each matrix as well as the informed step size $s_0^*$ given by the estimator in \Cref{fig:suitesparse}. In many cases, the adapted step size is $O(100)$ without relying on any preconditioning or matrix equilibration techniques. With a very large starting step size, the estimator with a threshold of $\Omega_\text{est}=10^7$ is also shown to be capable of adjusting the initial step size close to the adapted step size before the actual computation begins.

\begin{table}[htbp]
\centering
\caption{Properties of test matrices. The matrices are selected from a wide range of problem domains and each one is of size $N>10^6$.}
\begin{tabular}{c>{\centering\arraybackslash}p{0.4\textwidth}cc}
\toprule
\textbf{Matrix}       & \textbf{Problem Domain}                             & \textbf{$N$} & \textbf{Nonzeros} \\
\midrule
atmosmodd             & Computational Fluid Dynamics              & 1,270,432     & 8,814,880        \\
atmosmodj             & Computational Fluid Dynamics              & 1,270,432     & 8,814,880        \\
atmosmodl             & Computational Fluid Dynamics              & 1,489,752     & 10,319,760       \\
cage14                & Directed Weighted Graph                   & 1,505,785     & 27,130,349       \\
webbase-1M            & Directed Weighted Graph                   & 5,154,859     & 99,199,551       \\
G3\_circuit           & Circuit Simulation Problem                & 1,585,478     & 7,660,826        \\
Hamrle3               & Circuit Simulation Problem                & 1,447,360     & 5,514,242        \\
ecology1              & 2D/3D Problem                              & 1,000,000     & 4,996,000        \\
af\_shell10           & Structural Problem                        & 1,508,065     & 52,259,885       \\
Transport             & Structural Problem                        & 1,602,111     & 23,487,281       \\
CurlCurl\_3           & Model Reduction Problem                   & 1,219,574     & 13,544,618       \\
CurlCurl\_4           & Model Reduction Problem                   & 2,380,515     & 26,515,867       \\
nlpkkt80              & Optimization Problem                      & 1,062,400     & 28,192,672       \\
nlpkkt240             & Optimization Problem                      & 27,993,600    & 760,648,352      \\
thermal2              & Thermal Problem                           & 1,228,045     & 8,580,313        \\
\bottomrule
\end{tabular}
\label{tab:matrix_properties}
\end{table}

\section{Results on parallel scaling performance}\label{sec:scalability}

Numerical experiments in the previous section are conducted sequentially. This section investigates the parallel scalability of the proposed algorithm in a distributed-memory setting. 

We implement the adaptive s-step GMRES algorithm in \Cref{alg:adaptsgmres} with the scaled Newton polynomial basis. Additionally, we incorporate three additional variants of the GMRES algorithm to serve as benchmarks against our proposed methodology.
\begin{itemize}
  \item MGS-GMRES: a standard GMRES algorithm with Modified Gram-Schmidt (MGS) as the orthgonalization scheme. It is commonly observed that this variant exhibits suboptimal parallel performance, primarily attributed to the considerable communication costs incurred by the numerous global reductions involved.
  \item CGS2-GMRES: a GMRES variant that uses Classical Gram-Schmidt with re-orthogonalization (CGS2) as the orthogonalization scheme. Despite the nearly doubled computational cost compared to MGS, CGS2-GMRES is recognized for its numerical stability, attributed to the re-orthogonalization process~\cite{swirydowicz2020low}. Notably, this variant has $O(1)$ number of global reductions per iteration, often resulting in enhanced scalability compared to its MGS-based counterpart.
  \item BCGS2-TSQR-GMRES: the s-step GMRES algorithm that uses Block Classical Gram-Schmidt with reorthogonalization (BCGS2) as inter-ortho\-go\-na\-li\-zation and TSQR~\cite{demmel2012communication} as intra-ortho\-go\-na\-li\-zation. This variant closely resembles \Cref{alg:adaptsgmres} except that the intra-ortho\-go\-na\-li\-zation scheme is replaced with TSQR. TSQR is a communication-avoiding QR factorization scheme and is unconditionally stable. Given its superior numerical stability compared to CholQR, it is conjectured that this variant of GMRES maintains numerical stability, provided there exists an mechanism to actively monitor the condition number of Krylov basis matrices during intra-ortho\-go\-na\-li\-zation, such as via the incremental condition estimator proposed in \Cref{subsec:stopping}. While an in-depth numerical stability analysis of this GMRES variant is beyond the scope of this paper, our focus is primarily on its parallel performance for comparative purposes. To ensure fair comparisons, BCGS2-TSQR-GMRES employs the same scaled Newton polynomial basis as the adaptive s-step GMRES. 
\end{itemize}

Overall, we employ a total of four algorithms, comprising two column-wise GMRES variants (MGS-GMRES and CGS2-GMRES) and two s-step GMRES variants (BCGS2-TSQR-GMRES and Adaptive s-step GMRES). These algorithms are implemented in \texttt{C++} and utilize the \texttt{PETSc} library \cite{petsc-web-page} for parallel matrix/vector computations. The BLAS/LAPACK backend is provided by \texttt{Intel MKL}. Inter-process communication for Sparse Matrix-Vector (SpMV) multiplication and global reductions across the distributed network is facilitated by the built-in Message Passing Interface (MPI) support within \texttt{PETSc}. Given that the performance of TSQR heavily relies on its parallel implementation, we utilize the publicly available TSQR implementation in the open-source library \texttt{SLEPc} \cite{Hernandez:2005:SSF} for the sake of reproducibility. In \texttt{SLEPc}, the local QR factorizations in TSQR are carried out by LAPACK routines \texttt{geqrf} and \texttt{orgqr}. In addition, the Q matrix is computed explicitly in the \texttt{SLEPc} implementation of the TSQR kernel.

The codebase, including the libraries \texttt{PETSc 3.20} and \texttt{SLEPc 3.20}, is compiled using the \texttt{Intel oneAPI} compiler version \texttt{23.1.0}. Additionally, \texttt{Intel MPI} library version \texttt{21.9.0} is selected for parallel communication. The tests are performed on a cluster featuring Intel Xeon(R) Platinum 8280 CPUs with clock speed at 2.70GHz. Each node on the cluster is equipped with 56 cores and interconnected via Mellanox Infiniband HDR-100. We assign 1 MPI rank to each core as the mode of parallelization.

\subsection{Weak scaling}\label{sec:weakscaling}

We present weak scaling of our proposed algorithm using 3D Laplace matrices arising from 7-stencil finite difference discretization of uniform cubic grids. This matrix choice gives us the freedom to vary the dimension of the matrix without changing the sparsity pattern and its eigenvalue spectrum, allowing for fair comparisons. We vary the dimension of the matrix from $8 \times 10^6$ to $16\times 10^9$ and the number of nodes from 1 to 2048 (56 to 114,688 cores). For all benchmarks, we run 100 iterations without any restart. 

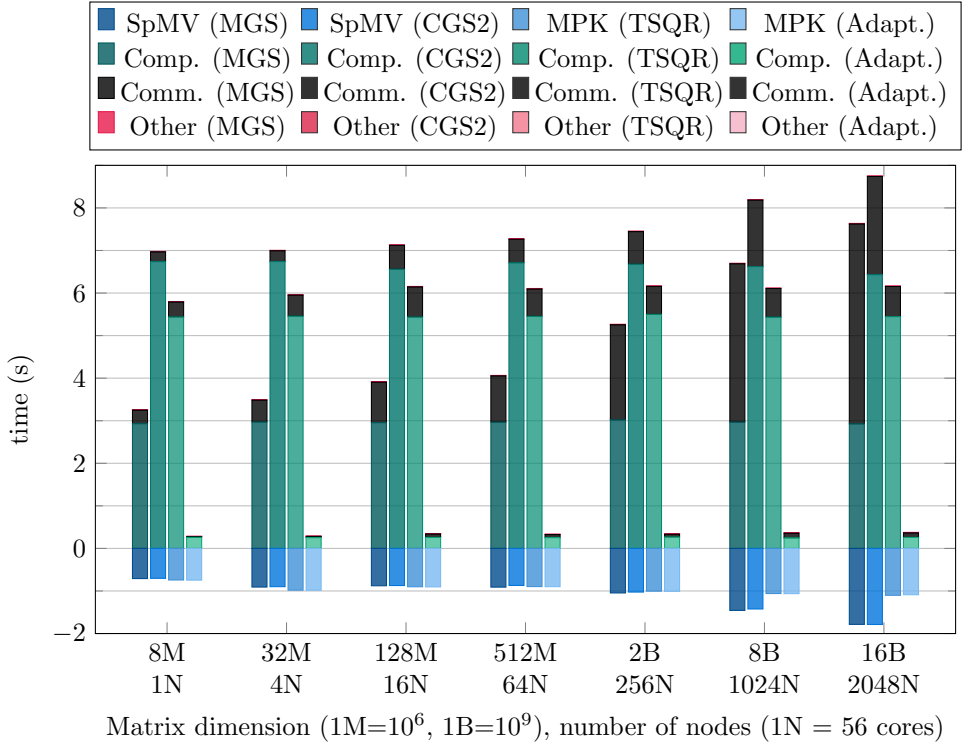
\begin{figure}
    \centering
    \begin{tikzpicture}
    \begin{axis}[
        xmode=log,
        width=\textwidth,
        height=.6\textwidth,
        bar width=.2cm,
        stack negative=separate,
        ybar stacked,
        ymin=-2,
        ymax=9,
        xtick={1,4,16,64,256,1024,4096},
        xticklabels={8M{\\}1N,32M{\\}4N,128M{\\}16N,512M{\\}64N,2B{\\}256N,8B{\\}1024N,16B{\\}2048N},
        xticklabel style={align=center,text width=10mm},
        xlabel={Matrix dimension (1M=$10^6$, 1B=$10^9$), number of nodes (1N = 56 cores)},
        x label style={at={(axis description cs:0.5,-0.15)},anchor=north},
        ylabel={time (s)},
        ymajorgrids=true,  
        yminorgrids=true,   
        minor y tick num=1,
        legend pos=outer north east,
        legend columns=4,
        legend style={at={(0.5,1.35)}, anchor=north},
    ]
    \addplot +[draw=black,
    color=blue2_2,
    fill opacity=0.8,
    bar shift=-.36cm] table [x=Node,y=SpMV_MGS] {Review/weak_L8.dat};
    \addlegendimage{draw=none,fill=blue3_2,opacity=0.8}; 
    \addlegendimage{draw=none,fill=blue4_2,opacity=0.8}; 
    \addlegendimage{draw=none,fill=blue5_2,opacity=0.8}; 
    \addplot +[draw=black,
    color=green2_2,
    fill opacity=0.8,
    bar shift=-.36cm] table [x=Node,y=Ortho_MGS] {Review/weak_L8.dat};
    \addlegendimage{draw=none,fill=green3_2,opacity=0.8}; 
    \addlegendimage{draw=none,fill=green4_2,opacity=0.8};
    \addlegendimage{draw=none,fill=green5_2,opacity=0.8};
    \addplot +[draw=black,
    color=black,
    fill opacity=0.8,
    bar shift=-.36cm] table [x=Node,y=Comm_MGS] {Review/weak_L8.dat};
    \addlegendimage{draw=none,fill=black,opacity=0.8}; 
    \addlegendimage{draw=none,fill=black,opacity=0.8};
    \addlegendimage{draw=none,fill=black,opacity=0.8};
    \addplot +[draw=black,
    color=hcred,
    fill opacity=0.8,
    bar shift=-.36cm] table [x=Node,y=Other_MGS] {Review/weak_L8.dat};
    % \addlegendimage{draw=none,fill=red2_1,opacity=0.8}; 
    \addlegendimage{draw=none,fill=red3_1,opacity=0.8};
    \addlegendimage{draw=none,fill=red4_1,opacity=0.8};
    \addlegendimage{draw=none,fill=red5_1,opacity=0.8};
    \legend{SpMV (MGS), SpMV (CGS2), MPK (TSQR), MPK (Adapt.), Comp. (MGS), Comp. (CGS2), Comp. (TSQR), Comp. (Adapt.), Comm. (MGS), Comm. (CGS2), Comm. (TSQR), Comm. (Adapt.), Other (MGS), Other (CGS2), Other (TSQR), Other (Adapt.)};
    \end{axis}
    \begin{axis}[
        xmode=log,
        width=\textwidth,
        height=.6\textwidth,
        bar width=.2cm,
        stack negative=separate,
        ybar stacked,
        ymin=-2,
        ymax=9,
        hide axis,
    ]
    \addplot +[draw=black,
    color=blue3_2,
    fill opacity=0.8,
    bar shift=-.12cm] table [x=Node,y=SpMV_CGS2] {Review/weak_L8.dat};
    \addplot +[draw=black,
    color=green3_2,
    fill opacity=0.8,
    bar shift=-.12cm] table [x=Node,y=Ortho_CGS2] {Review/weak_L8.dat};
    \addplot +[draw=black,
    color=black,
    fill opacity=0.8,
    bar shift=-.12cm] table [x=Node,y=Comm_CGS2] {Review/weak_L8.dat};
    \addplot +[draw=black,
    color=hcred,
    fill opacity=0.8,
    bar shift=-.12cm] table [x=Node,y=Other_CGS2] {Review/weak_L8.dat};
    \end{axis}
    \begin{axis}[
        xmode=log,
        width=\textwidth,
        height=.6\textwidth,
        bar width=.2cm,
        stack negative=separate,
        ybar stacked,
        ymin=-2,
        ymax=9,
        hide axis,
    ]
    \addplot +[draw=black,
    color=blue4_2,
    fill opacity=0.8,
    bar shift=.12cm] table [x=Node,y=SpMV_TSQR] {Review/weak_L8.dat};
    \addplot +[draw=black,
    color=green4_2,
    fill opacity=0.7,
    bar shift=.12cm] table [x=Node,y=Ortho_TSQR] {Review/weak_L8.dat};
    \addplot +[draw=black,
    color=black,
    fill opacity=0.8,
    bar shift=.12cm] table [x=Node,y=Comm_TSQR] {Review/weak_L8.dat};
    \addplot +[draw=black,
    color=hcred,
    fill opacity=0.8,
    bar shift=.12cm] table [x=Node,y=Other_TSQR] {Review/weak_L8.dat};
    \end{axis}
        \begin{axis}[
        xmode=log,
        width=\textwidth,
        height=.6\textwidth,
        bar width=.2cm,
        stack negative=separate,
        ybar stacked,
        ymin=-2,
        ymax=9,
        hide axis,
    ]
    \addplot +[draw=black,
    color=blue5_2,
    fill opacity=0.8,
    bar shift=.36cm] table [x=Node,y=SpMV_CholQR] {Review/weak_L8.dat};
    \addplot +[draw=black,
    color=green5_2,
    fill opacity=0.6,
    bar shift=.36cm] table [x=Node,y=Ortho_CholQR] {Review/weak_L8.dat};
    \addplot +[draw=black,
    color=black,
    fill opacity=0.8,
    bar shift=.36cm] table [x=Node,y=Comm_CholQR] {Review/weak_L8.dat};
    \addplot +[draw=black,
    color=hcred,
    fill opacity=0.8,
    bar shift=.36cm] table [x=Node,y=Other_CholQR] {Review/weak_L8.dat};
    \end{axis}
    \end{tikzpicture}
    \caption{Timing breakdown for weak scaling test of 3D Laplace matrix. Every four columns from left to right are: MGS-GMRES (MGS), CGS2-GMRES (CGS2), BCGS2-TSQR-GMRES (TSQR), adaptive s-step GMRES (Adapt.). The timing of SpMV/MPK is inverted and plotted below the x-axis in the bar plot to visually showcase the scalability of different operations without stacking all of them. The orthogonalization steps include both computational costs (Comp.) and communication costs (Comm.). Other miscellaneous operations such as updating the upper Hessenberg matrix and checking for convergence are included as well though their timings are almost negligible.}
    \label{fig:weakscaling}
\end{figure}

We first present a set of weak scaling experiments where no preconditioners are used so that all iterations are completed without premature convergence. For adaptive $s$-step GMRES algorithms, we use a step size of $s_0 = 100$ with scaled Newton polynomials. The initial step size estimator predicts $s_0^* = 100$ and the actual implementation produces the same step size as well in all cases. Hence, there is only one block iteration involved in each run of the adaptive $s$-step GMRES algorithm. The timing breakdown for SpMV and orthogonalization steps is shown in \Cref{fig:weakscaling}.

First, it is observed that the scalability of SpMV/MPK kernels, represented by the blue-themed bars in~\Cref{fig:weakscaling}, demonstrates consistent behavior across all parallel experiments. There is no discernible difference observed between the performance of SpMV in column-wise GMRES algorithms and that of MPK in s-step GMRES algorithms for most sizes. At a very large scale (1024 nodes or more), the MPK kernels in $s$-step GMRES algorithms seem to be slightly faster than SpMVs in column-based GMRES variants. Since MPK kernels have a similar number of floating point operations as SpMVs, the main difference comes from the fact that multiple SpMVs are called consecutively in MPKs before orthogonalization. Hence, we conjecture that the minor speedup is a result of better cache usage during the consecutive application of SpMVs in MPKs. Nevertheless, the SpMV/MPK kernels scale fairly well from 1 to 2048 nodes.

Furthermore, miscellaneous operations, including updating the upper Hessenberg matrix and checking for stopping criteria, are incorporated into the algorithm. However, their execution times are negligible compared to computationally intensive tasks such as SpMV multiplications and orthogonalizations. Consequently, the timings of these miscellaneous operations have minimal impact on the scaling performance of the algorithm.

Given the favorable scalability demonstrated by the SpMV/MPK kernel and minimal computational cost associated with other miscellaneous operations, our attention is directed towards the orthogonalization process for further discussion.

We compare the orthogonalization step among the four algorithms in~\Cref{fig:weakscaling} (green-themed bars). We first note that the orthogonalization step of MGS-GMRES demonstrates poor scalability as the number of nodes increases. This scalability issue primarily arises from the polynomial increase in the communication costs associated with a larger number of global reductions. It is worth recalling that for $s$ iterations, the number of global reductions in MGS grows quadratically as $O(s^2)$.

In contrast, the orthogonalization process in CGS2-GMRES exhibits improved scalability. This is attributed to the fact that the total number of global reductions in CGS2-GMRES is on the order of $O(s)$, resulting in a slower rate of increase in orthogonalization time with respect to the number of iterations $s$. However, it is important to note that the overall orthogonalization time in CGS2-GMRES is longer than that in MGS-GMRES. This is primarily due to the re-orthogonalization step in CGS2, which effectively doubles the computational cost of the CGS orthogonalization process. Consequently, the overall duration of the orthogonalization step in CGS2-GMRES is approximately twice as long as that in MGS-GMRES on a single node. The difference between the two diminishes with increasing number of nodes due to the overwhelming communication costs in MGS.

The performance of BCGS2-TSQR-GMRES surpasses that of CGS2-GMRES. BCGS2-TSQR demonstrates superior performance even on a single node compared to CGS2. This is largely attributed to speedup from the local QR LAPACK routines, \texttt{geqrf} and {orgqr}, in TSQR. Additionally, BCGS2-TSQR exhibits close-to-linear scalability by requiring only $O(1)$ global reductions for computing $s$ Krylov basis vectors. TSQR employs custom tree reductions of local QR factorizations. The theoretical communication costs are comparable to those of conventional global reductions. In practice, we do observe fairly constant communication costs across the weak scaling tests. However, it is worth highlighting that despite the enhanced performance of BCGS2-TSQR-GMRES compared to CGS2-GMRES, BCGS2-TSQR still falls short of MGS-GMRES for at one node. This discrepancy arises from the fact that BCGS2-TSQR involves approximately twice the number of floating-point operations as MGS-GMRES due to the re-orthogonalization process.

The adaptive s-step GMRES algorithm demonstrates the most efficient performance in the orthogonalization step among all algorithms. Due to its requirement of only $O(1)$ global reductions for $s$ basis vectors, the scalability of the algorithm approaches linearity, with the execution time remaining constant across 1 to 2048 nodes. One remarkable observation is that the orthogonalization operation is notably faster in the adaptive s-step GMRES compared to all other variants of GMRES. Comparing with BCGS2-TSQR-GMRES, our proposed algorithm primarily differs in the intra-orthogonalization scheme, where a (partial) CholQR method is employed instead of TSQR. Although CholQR and TSQR involve approximately the same number of floating-point operations~\cite{demmel2012communication}, CholQR exhibits better performance in practice. We conjecture that this is because CholQR predominantly utilizes BLAS-3 kernels, such as \texttt{gemm} and \texttt{trsm} routines, which are aggressively optimized for data locality and movement among various levels of cache memory. This superior performance of CholQR has also been observed in similar studies~\cite{yamazaki2015mixedchol}.

Overall, we demonstrate near-linear weak scaling of the adaptive $s$-step GMRES algorithm up to 2048 nodes (114,688 cores) for 3D Laplace matrices. In particular, the orthogonalization step in our proposed algorithm is capable of limiting the number of global reductions to $O(1)$ at the inter-node level and improving data locality at the intra-node level. 

\subsection{Additional weak scaling tests with preconditioners}
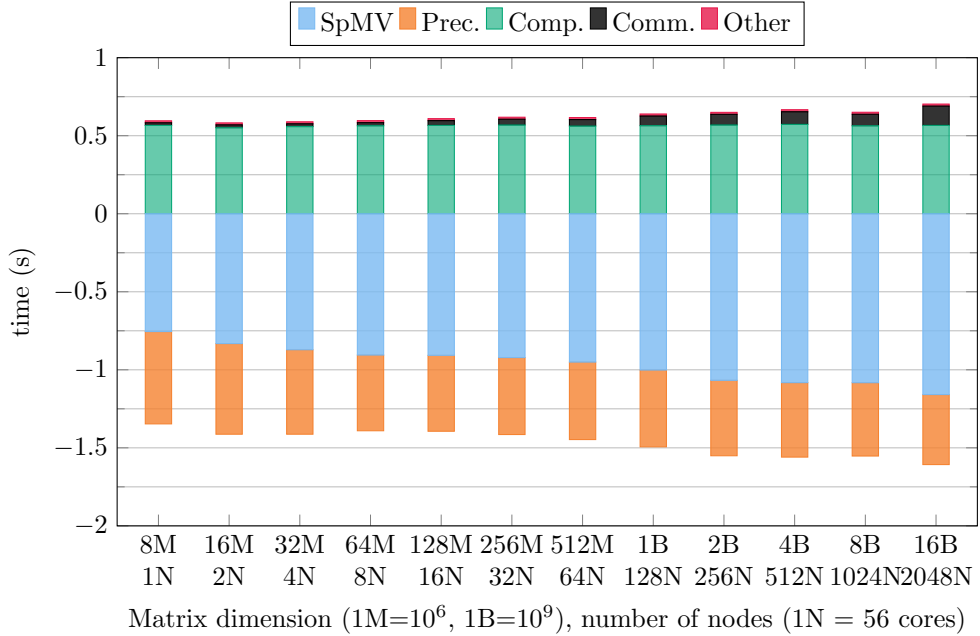
\begin{figure}
    \centering
    \begin{tikzpicture}
    \begin{axis}[
        xmode=log,
        width=\textwidth,
        height=.6\textwidth,
        bar width=.35cm,
        stack negative=separate,
        ybar stacked,
        xmin=0.667,
        xmax=3072,
        ymin=-2,
        ymax=1.0,
        xtick={1,2,4,8,16,32,64,128,256,512,1024,2048},
        xticklabels={8M{\\}1N,16M{\\}2N,32M{\\}4N,64M{\\}8N,128M{\\}16N,256M{\\}32N,512M{\\}64N,1B{\\}128N,2B{\\}256N, 4B{\\}512N, 8B{\\}1024N, 16B{\\}2048N},
        xticklabel style={align=center,text width=9.49mm},
        xlabel={Matrix dimension (1M=$10^6$, 1B=$10^9$), number of nodes (1N = 56 cores)},
        x label style={at={(axis description cs:0.5,-0.15)},anchor=north},
        ylabel={time (s)},
        ymajorgrids=true,  
        yminorgrids=true,
        minor y tick num=1,
        % legend pos=north east,
        legend columns=5,
        legend style={at={(0.5,1.12)}, anchor=north},
    ]
    \addplot +[draw=black,
    color=blue5_2,
    fill opacity=0.8,
    ] table [x=Node,y=SpMV] {Review/weak_L8_prec.dat};
    \addplot +[draw=black,
    color=hcorange,
    fill opacity=0.8,
    ] table [x=Node,y=Solve] {Review/weak_L8_prec.dat};
    \addplot +[draw=black,
    color=green5_2,
    fill opacity=0.6,
    ] table [x=Node,y=Ortho_Comp] {Review/weak_L8_prec.dat};
    \addplot +[draw=black,
    color=black,
    fill opacity=0.8,
    ] table [x=Node,y=Ortho_Comm] {Review/weak_L8_prec.dat};
    \addplot +[draw=black,
    color=hcred,
    fill opacity=0.8,
    ] table [x=Node,y=Others] {Review/weak_L8_prec.dat};
    \legend{SpMV, Prec., Comp., Comm., Other};
    \end{axis}
    \end{tikzpicture}
    \caption{Timing breakdown for weak scaling test of preconditioned 3D Laplace matrix using the adaptive s-step GMRES algorithm. The timings of SpMVs and preconditioning steps (Prec.) are inverted in the bar plot to showcase scalability of different operations. The orthogonalization steps include both computational costs (Comp.) and communication costs (Comm.). Other miscellaneous operations such as updating the upper Hessenberg matrix and checking for convergence are included as well though their timings are almost negligible.}
    \label{fig:weakscaling_pc}
\end{figure}

Here, we present an additional series of weak scaling experiments where we solve the identical 3D Laplace problem utilizing a local ILU(0) preconditioner. We run the same 100 iterations without any restart. In practice, the initial step size $s_0$ can reach 100 for this particular set of preconditioned numerical tests, but we purposely constrain the initial step size to $s_0=25$. This decision is made to demonstrate the scalability of our algorithm in situations where the adapted step size is not the same as the total number of iterations, necessitating multiple inter-orthogonalization steps. Overall, these tests are designed to emulate a more demanding numerical setting where a preconditioner is used and adapted step sizes are limited. Throughout all tests, the block size $s$ remains fixed at $s_0=25$ across iterations, with a total of four block iterations being executed. The local ILU(0) preconditioner is implemented using the built-in ILU preconditioner in PETSc. The results are summarized in \Cref{fig:weakscaling_pc}.

Our primary focus lies on the adaptive $s$-step GMRES algorithm, instead of other variants of GMRES. As depicted in~\Cref{fig:weakscaling_pc}, the costs associated with SpMVs remain comparable to those illustrated in~\Cref{fig:weakscaling}, since there are no modifications to the SpMV operations. In this series of scaling experiments, the preconditioning step (represented by orange-themed bars) exhibits relative constant runtime, attributable to the utilization of a local ILU(0) preconditioner that does not involve any communication. 

Our emphasis is therefore placed on the scalability of the orthogonalization step. With an increase in the number of nodes from 1 to 2048, the computational time of orthogonalization remains consistent, while communication costs see a gradual rise. Comparing the results against those of the adaptive $s$-step GMRES in~\Cref{fig:weakscaling}, a slight increment in the computational costs of orthogonalization is observed, rising from approximately $0.3s$ to $0.6s$. This is attributed to the reduction in step size from $s=100$ to $25$, resulting in a diminished utilization of cache. Nevertheless, communication costs remain minimal, at approximately $0.12s$ for 2048 nodes. The overall time for the orthogonalization process remains below $1s$. This still signifies a notable speedup compared to other variants of GMRES in~\Cref{fig:weakscaling}.

Overall, under the action of a preconditioner and with a reduced step size, the adaptive $s$-step GMRES algorithm still exhibits significant speedup compared to various GMRES variants and near-to-linear scaling performance up to 2048 nodes (114,688 cores).

\subsection{Strong scaling}
For strong scaling, we use a 3D Laplace matrix with a dimension of 64M on 1 to 64 nodes (56 to 3,584 cores). Similar to the weak scaling tests, we allow all benchmarks to run 100 iterations without any restart. A step size of 100 is used in all adaptive $s$-step GMRES algorithms. The initial step size estimator predicts $s_0^* = 100$ and agrees with the actual initial step size in all benchmarks. The time breakdown is shown in \Cref{fig:strongscaling}. Due to the large variation in magnitude, the corresponding speedup in log-scale is shown in \Cref{fig:strongscaling_speedup}.

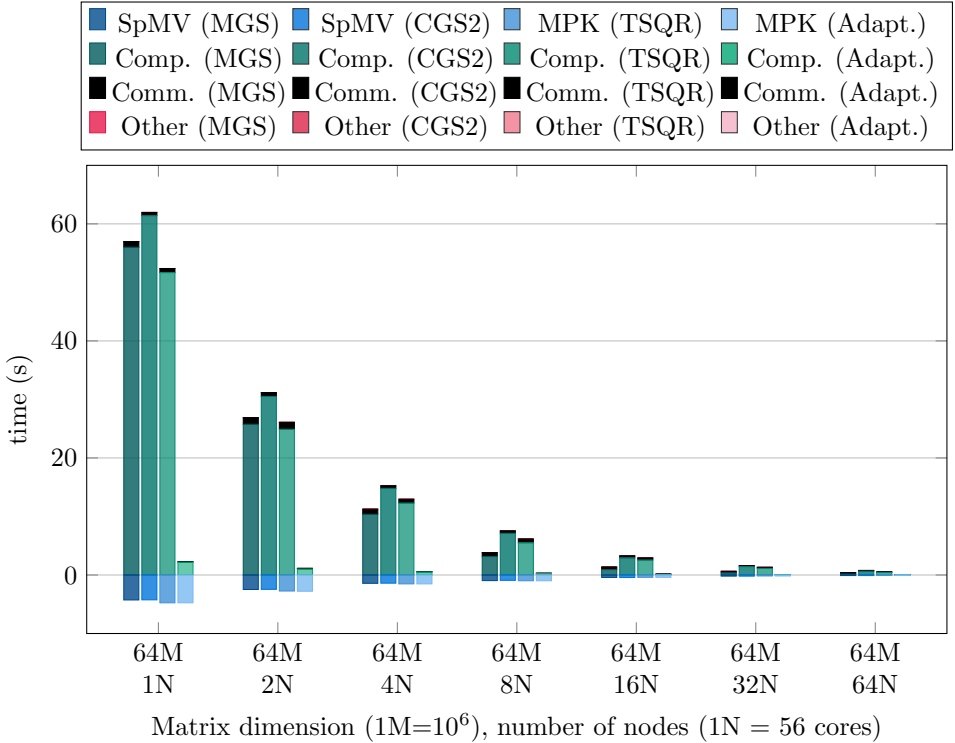
\begin{figure}[htbp]
    \centering
    \begin{tikzpicture}
    \begin{axis}[
        xmode=log,
        width=\textwidth,
        height=.6\textwidth,
        stack negative=separate,
        bar width=0.2cm,
        ybar stacked,
        ymin=-10,
        ymax=70,
        xtick={1,2,4,8,16,32,64},
        xticklabels={64M{\\}1N,64M{\\}2N,64M{\\}4N,64M{\\}8N,64M{\\}16N,64M{\\}32N,64M{\\}64N},
        xticklabel style={align=center,text width=10mm},
        xlabel={Matrix dimension (1M=$10^6$), number of nodes (1N = 56 cores)},
        x label style={at={(axis description cs:0.5,-0.15)},anchor=north},
        ylabel={time (s)},
        ymajorgrids=true,
        legend columns=4,
        legend style={at={(0.5,1.35)}, anchor=north},
    ]
    \addplot +[draw=black,
    color=blue2_2,
    fill opacity=0.8,
    bar shift=-.36cm] table [x=Nodes,y=SpMV_MGS] {Review/strong.dat};
    \addlegendimage{draw=none,fill=blue3_2,opacity=0.8}; 
    \addlegendimage{draw=none,fill=blue4_2,opacity=0.8}; 
    \addlegendimage{draw=none,fill=blue5_2,opacity=0.8}; 
    \addplot +[draw=black,
    color=green2_2,
    fill opacity=0.8,
    bar shift=-.36cm] table [x=Nodes,y=MGS] {Review/strong.dat};
    \addlegendimage{draw=none,fill=green3_2,opacity=0.8}; 
    \addlegendimage{draw=none,fill=green4_2,opacity=0.8};
    \addlegendimage{draw=none,fill=green5_2,opacity=0.8};
    \addplot +[draw=black,
    color=black,
    bar shift=-.36cm] table [x=Nodes,y=MGS_Comm] {Review/strong.dat};
    \addlegendimage{draw=none,fill=black}; 
    \addlegendimage{draw=none,fill=black};
    \addlegendimage{draw=none,fill=black};
    \addplot +[draw=black,
    color=hcred,
    fill opacity=0.8,
    bar shift=-.36cm] table [x=Nodes,y=Others_MGS] {Review/strong.dat};
    \addlegendimage{draw=none,fill=red3_1,opacity=0.8};
    \addlegendimage{draw=none,fill=red4_1,opacity=0.8};
    \addlegendimage{draw=none,fill=red5_1,opacity=0.8};
    \legend{SpMV (MGS), SpMV (CGS2), MPK (TSQR), MPK (Adapt.), Comp. (MGS), Comp. (CGS2), Comp. (TSQR), Comp. (Adapt.), Comm. (MGS), Comm. (CGS2), Comm. (TSQR), Comm. (Adapt.), Other (MGS), Other (CGS2), Other (TSQR), Other (Adapt.)};
    \end{axis}
    \begin{axis}[
        xmode=log,
        width=\textwidth,
        height=.6\textwidth,
        stack negative=separate,
        bar width=0.2cm,
        ybar stacked,
        ymin=-10,
        ymax=70,
        hide axis,
    ]
    \addplot +[draw=black,
    color=blue3_2,
    fill opacity=0.8,
    bar shift=-.12cm] table [x=Nodes,y=SpMV_CGS2] {Review/strong.dat};
    \addplot +[draw=black,
    color=green3_2,
    fill opacity=0.8,
    bar shift=-.12cm] table [x=Nodes,y=CGS2] {Review/strong.dat};
    \addplot +[draw=black,
    color=black,
    bar shift=-.12cm] table [x=Nodes,y=CGS2_Comm] {Review/strong.dat};
    \addplot +[draw=black,
    color=hcred,
    fill opacity=0.8,
    bar shift=-.12cm] table [x=Nodes,y=Others_CGS2] {Review/strong.dat};
    \end{axis}
    \begin{axis}[
        xmode=log,
        width=\textwidth,
        height=.6\textwidth,
        stack negative=separate,
        bar width=0.2cm,
        ybar stacked,
        ymin=-10,
        ymax=70,
        hide axis,
    ]
    \addplot +[draw=black,
    color=blue4_2,
    fill opacity=0.8,
    bar shift=.12cm] table [x=Nodes,y=SpMV_TSQR] {Review/strong.dat};
    \addplot +[draw=black,
    color=green4_2,
    fill opacity=0.7,
    bar shift=.12cm] table [x=Nodes,y=TSQR] {Review/strong.dat};
    \addplot +[draw=black,
    color=black,
    bar shift=.12cm] table [x=Nodes,y=TSQR_Comm] {Review/strong.dat};
    \addplot +[draw=black,
    color=hcred,
    fill opacity=0.7,
    bar shift=.12cm] table [x=Nodes,y=Others_TSQR] {Review/strong.dat};
    \end{axis}
    \begin{axis}[
        xmode=log,
        width=\textwidth,
        height=.6\textwidth,
        stack negative=separate,
        bar width=0.2cm,
        ybar stacked,
        ymin=-10,
        ymax=70,
        hide axis,
    ]
    \addplot +[draw=black,
    color=blue5_2,
    fill opacity=0.8,
    bar shift=.36cm] table [x=Nodes,y=SpMV_CholQR] {Review/strong.dat};
    \addplot +[draw=black,
    color=green5_2,
    fill opacity=0.6,
    bar shift=.36cm] table [x=Nodes, y=CholQR] {Review/strong.dat};
    \addplot +[draw=black,
    color=black,
    bar shift=.36cm] table [x=Nodes, y=CholQR_Comm] {Review/strong.dat};
    \addplot +[draw=black,
    color=hcred,
    fill opacity=0.6,
    bar shift=.36cm] table [x=Nodes, y=Others_CholQR] {Review/strong.dat};
    \end{axis}
    \end{tikzpicture}
    \caption{Timing breakdown for strong scaling test of 3D Laplace matrix. The dimension of the matrix is 64M. The legend carries the same meaning as the one in \Cref{fig:weakscaling}. 
    }
    \label{fig:strongscaling}
\end{figure}
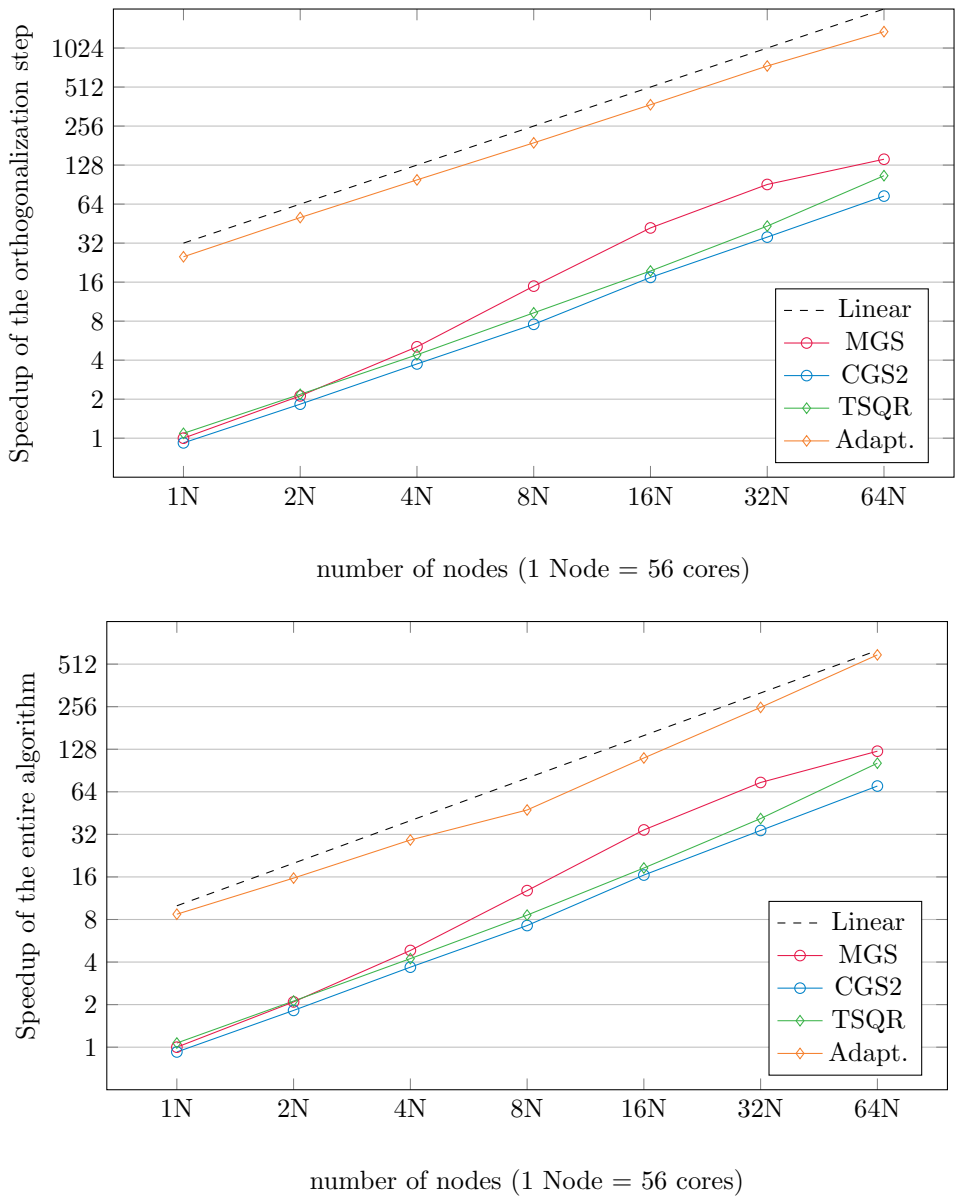
\begin{figure}[!htbp]
\centering
\subfloat{
    \begin{tikzpicture}
    \begin{axis}[
        xmode=log,
        ymode=log,
        width=.98\textwidth,
        height=.6\textwidth,
        stack negative=separate,
        ymin=0.5,
        ymax=2048,
        xtick={1,2,4,8,16,32,64},
        xticklabels={1N,2N,4N,8N,16N,32N,64N},
        xticklabel style={align=center,text width=10mm},
        xlabel={number of nodes (1 Node = 56 cores)},
        x label style={at={(axis description cs:0.5,-0.15)},anchor=north},
        ylabel={Speedup of the orthogonalization step},
        ymajorgrids=true,
        legend pos=north east,
        ytick={1,2,4,8,16,32,64,128,256,512,1024},
        yticklabels={1,2,4,8,16,32,64,128,256,512,1024},
        legend pos=south east,
    ]
    \addplot [dashed, color=black] coordinates {
    (1, 32)
    (64, 2048)
    };
    \addplot [mark=o, color=hcred] table [x=Nodes,y=MGS] {Review/strong_speedup.dat};
    \addplot [mark=o, color=hcblue] table [x=Nodes,y=CGS2] {Review/strong_speedup.dat};
    \addplot [mark=diamond,color=hcgreen] table [x=Nodes,y=TSQR] {Review/strong_speedup.dat};
    \addplot [mark=diamond,color=hcorange] table [x=Nodes,y=CholQR] {Review/strong_speedup.dat};
    \legend{Linear, MGS, CGS2, TSQR, Adapt.};
    \end{axis}
    \end{tikzpicture}
}\\
\subfloat{
    \begin{tikzpicture}
    \begin{axis}[
        xmode=log,
        ymode=log,
        width=.98\textwidth,
        height=.6\textwidth,
        stack negative=separate,
        ymin=0.5,
        ymax=1024,
        xtick={1,2,4,8,16,32,64},
        xticklabels={1N,2N,4N,8N,16N,32N,64N},
        xticklabel style={align=center,text width=10mm},
        xlabel={number of nodes (1 Node = 56 cores)},
        x label style={at={(axis description cs:0.5,-0.15)},anchor=north},
        ylabel={Speedup of the entire algorithm},
        ymajorgrids=true,
        legend pos=north east,
        ytick={1,2,4,8,16,32,64,128,256,512},
        yticklabels={1,2,4,8,16,32,64,128,256,512},
        legend pos=south east,
    ]
    \addplot [dashed, color=black] coordinates {
    (1, 10)
    (64, 640)
    };
    \addplot [mark=o, color=hcred] table [x=Nodes,y=Total_MGS] {Review/strong_speedup.dat};
    \addplot [mark=o, color=hcblue] table [x=Nodes,y=Total_CGS2] {Review/strong_speedup.dat};
    \addplot [mark=diamond,color=hcgreen] table [x=Nodes,y=Total_TSQR] {Review/strong_speedup.dat};
    \addplot [mark=diamond,color=hcorange] table [x=Nodes,y=Total_CholQR] {Review/strong_speedup.dat};
    \legend{Linear, MGS, CGS2, TSQR, Adapt.};
    \end{axis}
    \end{tikzpicture}
}
\caption{Speedup plot for strong scaling test of 3D Laplace matrix with dimension 64M. The four algorithms (and their corresponding legends) are MGS-GMRES (MGS), CGS2-GMRES (CGS2), BCGS2-TSQR-GMRES (TSQR) and adaptive s-step GMRES (Adapt.). Top figure shows the speedup of orthogonalization operations (including both computation and commmunication) and the bottom figure shows the speedup of overall runtime. Both algorithms use the corresponding MGS-GMRES time at 1 node as the reference for computing the speedup.}
\label{fig:strongscaling_speedup}
\end{figure}

The main observation here is that even at 1 node, the adaptive $s$-step GMRES algorithm significantly outperforms all other GMRES variants. Most of the speedup comes from the difference in computational cost of orthogonalization operations, as adaptive $s$-step GMRES algorithm makes use of BLAS-3 kernels with better data locality. This difference is consistently observed from 1 node to 64 nodes. 

The trend is better represented in \Cref{fig:strongscaling_speedup}. As shown in the top figure of \Cref{fig:strongscaling_speedup}, orthogonalization costs (both computation and communication) in standard MGS-GMRES shows linear scaling from 1 to 2 nodes at first, followed by superlinear scaling up to 16 nodes due to increased amount of cache memory, and then degraded performance at 64 nodes due to overwhelming communication costs. 
On the other hand, orthogonalization operations for the all other variants show close-to-linear speedup. We would like to remark that although both CGS2 and BCGS2-TSQR orthogonalization schemes have approximately twice the amount of floating point operations as MGS, their runtime is close to the runtime of MGS due to the use of BLAS-2 and BLAS-3 kernels that are better optimized for cache memory. Similar to the results in weak scaling tests, the orthogonalization scheme in adaptive s-step GMRES outperforms all other variants due to its use of highly optimized BLAS-3 kernels. An approximately $25\times$ improvement over MGS is observed at one node.

In addition, all algorithms have similar SpMV/MPK costs, as they compute the same number of SpMVs. In consequence, the scaling performance of the overall algorithm in the bottom figure of~\Cref{fig:strongscaling_speedup} closely follows the scaling performance of orthogonalization operations. In general, the adaptive $s$-step GMRES algorithm shows linear (strong) scaling up to 64 nodes with significant improvement over other variants of the GMRES algorithm.

\section{Conclusions}
\label{sec:conclusions}
We have developed a communication-avoiding adaptive $s$-step GMRES algorithm that is numerically stable. The algorithm uses a partial CholQR procedure that ensures that the block orthogonalization schemes are carried out within the stability limit by adjusting the step size on the fly. The overall orthogonalization produces minimal LOO errors, which are near the machine epsilon. With the introduction of scaled Newton polynomials, parallel scaling performance of the proposed algorithm in a distributed memory architecture is presented, showing near-linear scaling of the adaptive $s$-step GMRES algorithm up to 114,688 cores.

In addition to performance results, we also used the adapted step size as a metric to evaluate different variants of the algorithm. We investigated the stability of different polynomial bases for MPK and different preconditioning techniques. The robustness of scaled Newton polynomials is demonstrated, and the influence of the eigenvalue distribution on the adapted step size is emphasized. The initial step size estimator is developed based on scaled Newton polynomials to mitigate the amount of wasted computations and communications and has been shown to give well-informed step sizes close to the actual step size limit.

\section*{Acknowledgments}
The authors gratefully acknowledge the support of NASA Transformational Tools and Technologies program under NASA Grant 80N\-SSC\-18M\-0152, and contributions of William K. Anderson, project manager.
The authors also acknowledge the Texas Advanced Computing Center (TACC) at The University of Texas at Austin for providing HPC resources that have contributed to the research results reported within this paper.

\bibliographystyle{siamplain}
\bibliography{references}

\appendix

\end{document}